\documentclass[%
11pt,
onecolumn,
tightenlines,
superscriptaddress,
preprintnumbers,
nofootinbib,
amsmath,amssymb,amsthm,
physrev,
eqsecnum,
]{revtex4-2}

\usepackage{isomath}
\usepackage{amsmath,amsthm}
\usepackage{amsbsy}
\usepackage{amssymb}
\usepackage{amscd}
\usepackage{amsfonts}
\usepackage{stmaryrd}
\usepackage{euscript}
\usepackage[utf8]{inputenc}
\usepackage[T1]{fontenc}
\usepackage{newtxtext} 
\everymath{\displaystyle}
\usepackage{exscale}

\usepackage{graphicx}
\graphicspath{ {media/} }

\usepackage{hyperref}

\usepackage[small]{titlesec}

\titlespacing*{\section}{0pt}{12pt plus 4pt minus 2pt}{2pt plus 2pt minus 2pt}
\titlespacing*{\subsection}{0pt}{12pt plus 4pt minus 2pt}{2pt plus 2pt minus 2pt}
\titlespacing*\subsubsection{0pt}{12pt plus 4pt minus 2pt}{2pt plus 2pt minus 2pt}
\titlespacing*\paragraph{0pt}{12pt plus 4pt minus 2pt}{2pt plus 2pt minus 2pt}

\makeatletter

    \renewcommand*{\p@subsection}{}
    
    \renewcommand*{\p@subsubsection}{}
\makeatother

\usepackage{isomath}
\usepackage{amsmath}
\usepackage{amssymb}
\usepackage{amscd}
\usepackage{amsfonts}
\usepackage{xcolor}
\usepackage{soul}

\newtheorem{remark}{Remark}[section]

\newcommand{\bfchi}{\mathbold {\chi}}

\newcommand{\bfsigma}{\mathbold {\sigma}}

\newcommand{\bfxi}{\mathbold {\xi}}

\DeclareMathOperator{\divergence}{div}

\DeclareMathOperator{\trace}{tr}

\newcommand{\const}{\mathrm{const.}}

\newcommand{\parderiv}[2]{\frac{\partial #1}{\partial #2}}

\newcommand{\dm}{\ \mathrm{d}}

\newcommand{\deriv}[2]{\frac{\dm #1}{\dm #2}}

\newcommand{\bfb}{{\mathbold b}}

\newcommand{\bfe}{{\mathbold e}}

\newcommand{\bfg}{{\mathbold g}}

\newcommand{\bfn}{{\mathbold n}}

\newcommand{\bft}{{\mathbold t}}

\newcommand{\bfv}{{\mathbold v}}

\newcommand{\bfx}{{\mathbold x}}

\newcommand{\bfB}{{\mathbold B}}

\newcommand{\bfF}{{\mathbold F}}

\newcommand{\bfI}{{\mathbold I}}

\newcommand{\bfP}{{\mathbold P}}

\newcommand{\bfV}{{\mathbold V}}

\newcommand{\bfX}{{\mathbold X}}

\newcommand{\bfFr}{\bfF_{relax}}

\begin{document}

\preprint{To appear in Mathematics and Mechanics of Solids (DOI: \url{https://doi.org/10.1177/10812865211054573})}

\title{\Large{Accretion and Ablation in Deformable Solids with an Eulerian Description: Examples using the Method of Characteristics}}

\author{S. Kiana Naghibzadeh}
    \email{snaghibz@andrew.cmu.edu}
    \affiliation{Department of Civil and Environmental Engineering, Carnegie Mellon University}

\author{Noel Walkington}
    \email{noelw@cmu.edu}
    \affiliation{Center for Nonlinear Analysis, Department of Mathematical Sciences, Carnegie Mellon University}
    
\author{Kaushik Dayal}
    \email{Kaushik.Dayal@cmu.edu}
    \affiliation{Department of Civil and Environmental Engineering, Carnegie Mellon University}
    \affiliation{Center for Nonlinear Analysis, Department of Mathematical Sciences, Carnegie Mellon University}
    \affiliation{Department of Materials Science and Engineering, Carnegie Mellon University}
    
\date{\today}

\begin{abstract}
    Accretion and ablation, i.e. the addition and removal of mass at the surface, is important in a wide range of physical processes including solidification, growth of biological tissues, environmental processes, and additive manufacturing.
    The description of accretion requires the addition of new continuum particles to the body, and is therefore challenging for standard continuum formulations for solids that require a reference configuration.
    Recent work has proposed an Eulerian approach to this problem, enabling the side-stepping of the issue of constructing the reference configuration.
    However, this raises the complementary challenge of determining the stress response of the solid, which typically requires the deformation gradient that is not immediately available in the Eulerian formulation.
    To resolve this, the approach introduced the elastic deformation as an additional kinematic descriptor of the added material, and its evolution has been shown to be governed by a transport equation.
    In this work, the method of characteristics is applied to solve concrete simplified problems motivated by biomechanics and manufacturing.
    Specifically, (1) for a problem with both ablation and accretion in a fixed domain, and (2) for a problem with a time-varying domain, the closed-form solution is obtained in the Eulerian framework using the method of characteristics without explicit construction of the reference configuration.
\end{abstract}

\maketitle


\section{Introduction} \label{sec:intro}

Accretion or surface growth involves the addition of mass at the boundary of a body.
In a continuum setting, this requires the introduction of new material particles to the body \cite{skalak1982analytical}, and poses unusual challenges.
Surface growth occurs in a range of settings including in biological tissue \cite{taber1995biomechanics}, planetary formation \cite{brown1963gravitational, kadish2008stresses}, solidification and casting processes \cite{schwerdtfeger1998stress}, etching process in silicon wafers \cite{rao2000modelling}, additive manufacturing processes \cite{drozdov1998continuous}, construction of masonry structures \cite{bacigalupo2012effects}, and cell motility via assembly and disassembly of Actin networks \cite{papadopoulos2010surface}.

The addition and removal of material points in surface growth causes difficulties in the usual procedure of defining a fixed set of points as the reference configuration. 
Dealing with this in linear elasticity is somewhat easier \cite{ong2004equations, bacigalupo2012effects,naumov1994mechanics,kadish2005stresses,brown1963gravitational}.
However, in nonlinear descriptions of surface growth, the definition of the reference configuration is more subtle.
Pioneering work by Skalak, Hoger, and co-workers provided a systematic approach to study the kinematics of surface growth \cite{skalak1982analytical,skalak1997kinematics}. 
However, these studies did not consider deformation, and assumed that the body was rigid. 
Recent work has extended these ideas to account for deformation and stress \cite{sozio2017nonlinear, sozio2020nonlinear,sozio2019nonlinear,tomassetti2016steady,abi2019kinetics,abi2020surface,abeyaratne2020treadmilling,von2020morphogenesis,zurlo2017printing,truskinovsky2019nonlinear}.
These studies used a Lagrangian formulation, i.e., roughly the kinematic description is based on reference particles in an \textit{evolving} reference configuration.

To avoid the difficulties of dealing with an evolving reference configuration, our previous work \cite{naghibzadeh2021surface} developed an Eulerian description of surface growth.
This has the advantage that we do not need to explicitly calculate the reference configuration, e.g. \cite{clayton2013nonlinear,clayton2019nonlinear,clayton2020aplace}.
On the other hand, in a solid, the stress response requires knowledge of the deformation gradient $\bfF$, and the latter is not straightforward to calculate in an Eulerian setting.
Prior work, notably \cite{liu2001eulerian,kamrin2009eulerian,kamrin2012reference}, have developed approaches to find the deformation gradient in an Eulerian setting; specifically, \cite{liu2001eulerian} uses an evolution equation for $\bfF$ that does not require explicit calculation of the reference configuration.

While the Eulerian approach works well for the standard setting of a fixed set of material particles without growth, there is the challenge of defining the reference state and kinematic information of the added material particles during surface growth.
Specifically, the specification of $\bfF$ for the added material must ensure that kinematic compatibility is not violated.
To avoid dealing with the issue of kinematic compatibility, we introduced two additional kinematic descriptors: the relaxed zero-stress deformation $\bfF_{relax}$ and the elastic deformation $\bfF_e$.
These are related to the deformation gradient through the relation $\bfF_e = \bfF \bfF_{relax}^{-1}$.
They enable us to eliminate $\bfF$ -- as well as $\bfF_{relax}$, using that it is transported with the material particles -- and thereby formulate a model that is posed in terms of $\bfF_e$ as the sole kinematic descriptor for the stress response.

The evolution of $\bfF_e$ is governed by a transport equation, and transparently handles both accretion and ablation: accretion corresponds to inflow, and requires the specification of boundary conditions, while ablation corresponds to outflow and requires no boundary conditions.
Further, $\bfF_e$ has no requirements of kinematic compatibility, enabling a simple and direct approach to the specification of boundary conditions.
The boundary conditions are specified in terms of $\bfF_e$ for the added particles, which is directly and solely related to their given stress state.
We can think of $\bfF_e$ and $\bfFr$ as providing a mechanism for the added material to define its own reference configuration that satisfies kinematic compatibility.

We apply this formulation to examine two simplified concrete problems of recent interest.
In the first problem, motivated by the assembly and disassembly of actin micro-filaments during polymerization and depolymerization on a spherical bead \cite{pantaloni2001mechanism,noireaux2000growing}, and building on the model of this phenomenon from \cite{tomassetti2016steady}, we examine simultaneous accretion at the fixed inner surface and ablation at the time-dependent outer surface of a hollow sphere.
In the second problem, motivated by the wound roll process of additive manufacturing which has been applied to processes ranging from the packaging of flexible films in paper and magnetic tape industries \cite{yagoda1980resolution} to tissue engineered blood vessels \cite{konig2009mechanical}, we study the deformation and stress distribution when stress-free particles are deposited at the outer layer of a hollow cylinder with a time-dependent pressure applied at the inner surface.
We solve these problems in closed-form with the method of characteristics, using that the evolution of $\bfF_e$ is governed by a transport equation.

\paragraph*{Organization.}
Section \ref{sec:PhysicalModel} summarizes the approach proposed in \cite{naghibzadeh2021surface}.
We start with a discussion on balance laws and jump / boundary conditions, using the view of surface growth as a localized source on the growing surface.
We then turn to the kinematics of the deformation in the Eulerian setting with accretion.
Sections \ref{sec:inne_acc} and \ref{sec:outer_acc} solve specific examples of (1) simultaneous accretion at the fixed inner surface and ablation at the time-dependent outer surface of a hollow sphere, (2) accretion on the outer surface of a long hollow cylinder with a time-dependent pressure applied at the inner surface. We used the method of characteristics to solve the transport equation that governs the evolution of the elastic deformation.



\section{Physical model and governing equations}
\label{sec:PhysicalModel}

\paragraph*{Notation.}

The motion is denoted by $\bfx = \bfchi(\bfX , t)$, where $\bfX$ is the location in the reference configuration of a point that is currently at the spatial location $\bfx$.
Also, the inverse of the motion $\bfchi^{-1} (\bfx , t)$  maps the current spatial location of a particle at $\bfx$ to the location $\bfX$ in the reference configuration. 
The deformation gradient is $\bfF = \parderiv{\bfchi (\bfX , t)}{\bfX}$.
For short, we write $\bfF(\bfx,t) = \bfF(\bfchi^{-1}(\bfX,t),t)$.

Because most of our work is in the Eulerian setting, all differential operators (e.g. $\nabla, \divergence$) imply derivatives with respect to $\bfx$, except where explicitly stated.
Similarly, the argument of various field quantities will implicitly be $\bfx$, except where explicitly stated.

\subsection{Balance laws}

Following \cite{naghibzadeh2021surface}, we model surface growth through the introduction of sources that are localized on the growing surface.
The source terms can be understood as a coarse-grained approach to treat the complex process of growth without considering the fine details of the processes in the ambient environment outside the growing body; the growth is defined only in terms of net mass and linear momentum transfer.
To define the source terms, we specify the rate of mass addition per unit area, $M$; and the vectorial rate of linear momentum addition per unit area, $\bfP$.
We can then infer the velocity (or specific momentum) of the added material at the instant of attachment, $\bfv_a := \bfP / M$; and the mass density, $\rho := \frac{M}{\bfv_a \cdot \hat\bfn}$.
Accretion is modeled by $M>0$ and ablation by $M<0$.
We assume that the added particles do not carry angular momentum such as due to individual particle spins, i.e., there is no source of angular momentum due to growth.

The source terms appear in the jump conditions that represent the balances of mass and momentum at the growing surface:
\begin{align}
    \text{Mass: } & \llbracket \rho (\bfV_b - \bfv) \cdot \hat{\bfn} \rrbracket  = M
    \\
    \text{Momentum: } & \llbracket \rho \bfv ( (\bfV_b - \bfv) \cdot \hat{\bfn} )\rrbracket + \llbracket \bfsigma \hat{\bfn} \rrbracket = M \hat{\bfv}_a
\end{align}
We have used the jump operator $\llbracket \cdot \rrbracket$ to denote the difference between the limits when approaching the surface of discontinuity from either side; $\bfV_b$ is the velocity of the surface; and $\hat\bfn$ is the surface normal.

We make the assumption that the ambient environment outside of the growing body has a negligible effect on the mechanics of the body.
This leads to the significant simplification that the domain of the problem now involves only the solid body, and the jump conditions become boundary conditions:
\begin{align}
     \label{eqn:bc-balance-mass}
    \text{Mass: } & \rho (\bfV_b - \bfv) \cdot \hat{\bfn}  = M
    \quad \Longleftrightarrow \quad
    \bfV_b \cdot \hat{\bfn} = \bfv \cdot \hat{\bfn}  + \frac{M}{\rho}
    \\
     \label{eqn:bc-balance-momentum}
    \text{Momentum: } & \rho \bfv ((\bfV_b - \bfv) \cdot \hat{\bfn}) + 
    \bfsigma \hat{\bfn} - \bft_b = 
    M \bfv_a 
    \quad \Longleftrightarrow \quad
    \bfsigma \hat{\bfn} 
    =
    M (\bfv_a - \bfv) + \bft_b
\end{align}
where $\bft_b$ is the traction due to the external forces at the boundary of the growing body.
We further assume that the inertial terms $\rho \bfv ((\bfV_b - \bfv) \cdot \hat{\bfn})$ and $M \bfv_a$ are negligible compared to $\bfsigma$.
Then, \eqref{eqn:bc-balance-momentum} reduces to $\bfsigma \hat{\bfn} = \bft_b$.

\begin{figure}[htb!]
    \includegraphics[width=\textwidth]{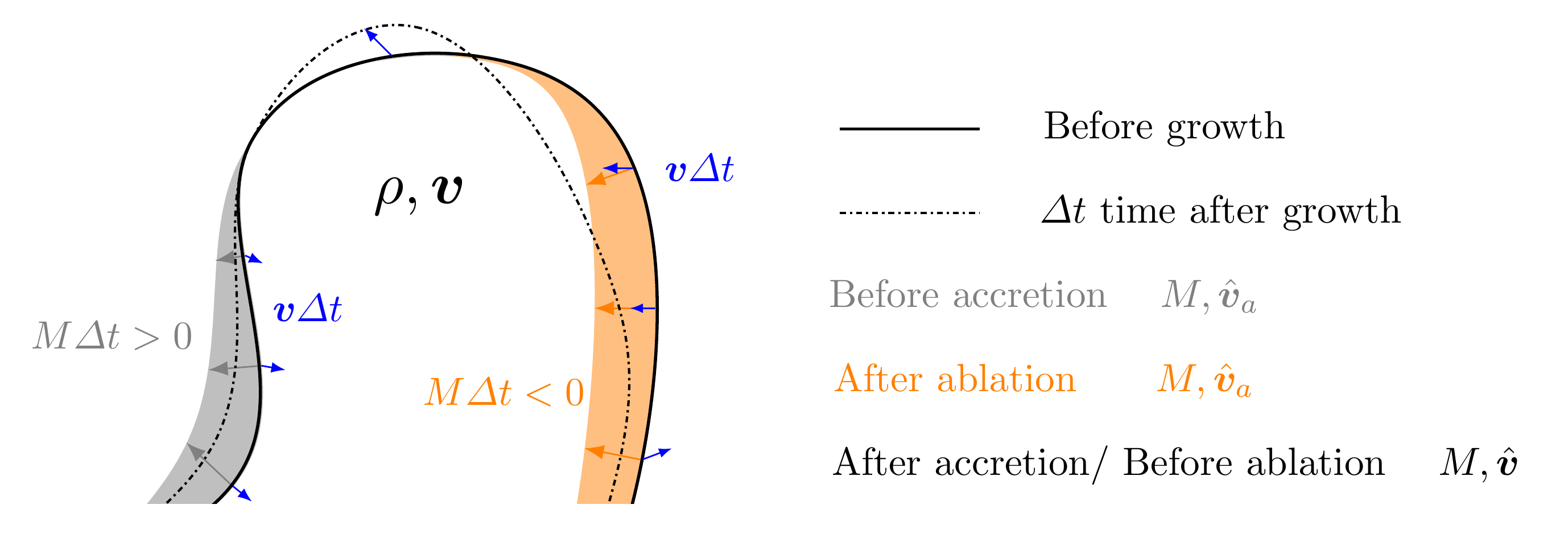}
    \caption{A schematic of the evolution of a growing body in a small interval of time $\Delta t$. The spatial location of the growing boundary depends on the growth velocity and the continuum particle velocity at the boundary.}
    \label{fig:definitions}
\end{figure}

The balance laws in the bulk are standard\footnote
{
    The balance of angular momentum provides only the standard result that the Cauchy stress must be symmetric. 
}:
\begin{align}
    \label{eqn:bulk-balance-mass}
    \text{Mass: } & \parderiv{\rho}{t} + \divergence(\rho \bfv) = 0
    \\
    \label{eqn:bulk-balance-momentum}
    \text{Momentum: } & \parderiv{\ }{t} (\rho \bfv) + \divergence (\rho \bfv \otimes \bfv)= \rho \bfb + \divergence(\bfsigma)
\end{align}
where $\bfsigma$ is the Cauchy stress; $\rho$ is the mass density; $\bfv$ is the particle velocity; $\bfb$ is the body force; $\bft_b$ is the traction vector due to the external forces at the boundary of the growing body; and $\bfV_b$ is the boundary velocity. From \eqref{eqn:bc-balance-mass}, the boundary velocity $\bfV_b$ has contributions from $M/\rho$ and $\bfv$ (Figure \ref{fig:definitions}).

\subsection{Kinematics}

The stress response of the material generally requires knowledge of the deformation gradient $\bfF = \parderiv{\bfchi (\bfX , t)}{\bfX} = \left( \parderiv{\bfchi^{-1} (\bfx , t)}{\bfx} \right)^{-1}$; we notice that it requires knowledge of the motion or its inverse.
However, we adapt the approach from \cite{liu2001eulerian} that formulates the evolution of $\bfF$ without requiring us to solve for the deformation map.

    We can write the material time derivative of $\bfF$ as:
    \begin{equation}
        \parderiv{\bfF(\bfX,t)}{t}
        =
        \parderiv{\ }{t}\parderiv{\bfchi(\bfX,t)}{\bfX}
        =
        \parderiv{\ }{\bfX}\parderiv{\bfx}{t}
        =
        \parderiv{\bfv}{\bfX}
        =
        \parderiv{\bfv}{\bfx}\parderiv{\bfx}{\bfX}
        =
        \left(\nabla\bfv\right) \bfF
        \Rightarrow
        \dot{\bfF}=\nabla\bfv \bfF
    \end{equation}
    where $\dot{\bfg}$ represents the material time derivative of a quantity $\bfg$.
It follows that the evolution of $\bfF$ satisfies the transport equation:
\begin{equation}
    \label{eqn:def-grad-transport}
    \parderiv{\bfF(\bfx,t)}{t}+(\bfv\cdot\nabla)\bfF=\left(\nabla\bfv \right) \bfF
\end{equation}
While this equation preserves the kinematic compatibility of $\bfF$ for a fixed set of material particles, the setting of surface growth requires imposing boundary conditions -- corresponding to the deformation state of the added material -- that must also be consistent with kinematic compatibility.
Assume that the stress state of the added particles is given to be $\bfsigma^*$.
Denote the stress response function of the added particles by $\hat\bfsigma(\cdot)$, and assume for simplicity and without loss of generality that $\hat\bfsigma(\bfI) = \bf 0$.
We notice immediately that $\bfsigma^* = \hat\bfsigma(\bfF)$; assuming for simplicity that the stress response is invertible\footnote{
    The stress response can be invertible only up to a rotation, but we can make any convenient choice for the unconstrained rotation.
} implies that $\bfF = \hat\bfsigma^{-1}(\bfsigma^*)$ is completely determined at the time of attachment.
Since $\bfsigma^*$ is specified independently of the state of the body, this could potentially violate kinematic compatibility at the time of attachment.

We therefore introduce two additional kinematic descriptors: $\bfF_{relax}$ that quantifies the relaxed (i.e., zero stress) shape of the added particles; and $\bfF_e := \bfF \bfFr^{-1}$ that quantifies the elastic deformation\footnote{
    These kinematic quantities are very similar to those proposed in  in surface growth by \cite{sozio2019nonlinear,zurlo2017printing,truskinovsky2019nonlinear}, and builds on ideas from the bulk setting that go back many decades \cite{sadik2017origins,garikipati2009kinematics}.
}.
We highlight that neither $\bfF_e$ nor $\bfF_{relax}$ need satisfy any requirements of kinematic compatibility.
We next redefine the stress response function to be $\hat\bfsigma(\bfg) \mapsto \hat\bfsigma(\bfg \bfF_{relax}^{-1})$.
We have the obvious interpretation of $\bfF_e$ as the part of $\bfF$ that causes stress, i.e. $\bfsigma = \hat\bfsigma(\bfF_e)$. 
This corresponds to redefining the reference of the added material.
Prior to attachment, the deformation state is irrelevant, and post-attachment deformation is governed by the standard equations of continuum mechanics.
For an unattached particle, the reference can be changed in arbitrary ways and, in particular, it can be changed to ensure that we satisfy kinematic compatibility.
In this regard, we highlight that the stress state $\bfsigma^*$ of the added material is a well-defined and physically-meaningful quantity, whereas the deformation can be chosen for convenience by appropriately defining the reference at the time of attachment.

    Substituting $\bfF = \bfF_e \bfFr$ in \eqref{eqn:def-grad-transport}, we find:
    \begin{align}
        \label{eqn:def-elastic-evolution}
        \parderiv{\bfF_e(\bfx,t)}{t}+(\bfv\cdot\nabla)\bfF_e=\left(\nabla\bfv \right) \bfF_e
    \end{align}
    We have assumed that $\bfFr$ is constant in time for each material particle, i.e., $\parderiv{\bfF_{relax}}{t} + (\bfv \cdot \nabla)\bfF_{relax}= 0$; \cite{naghibzadeh2021surface} discusses the more general case where $\bfFr$ can evolve in time.

To compute the stress response of the body, we only need to know $\bfF_e$, so we solve only \eqref{eqn:def-elastic-evolution}.
This is a transport equation, and requires boundary conditions only at inflow boundaries \cite{trangenstein2009numerical}, i.e., where the velocity $\bfv$ is inwards with respect to the \textit{moving} boundary.
Using the given stress state $\bfsigma^*$ of the added particles, we specify $\bfF_e = \hat\bfsigma^{-1}(\bfsigma^*)$ at the accretion boundaries.
    While $\bfF_e$ has 9 components that must be prescribed at the inflow boundary, we highlight that $\bfF_e$ is prescribed by the stress only up to a free rotation.
    Loosely, the stress provides a prescription of only 6 components, broadly similar to the approach by Truskinovsky \cite{truskinovsky2019nonlinear,zurlo2017printing} that projects stresses.

No boundary condition on $\bfF_e$ are needed for the ablation boundaries \cite{trangenstein2009numerical}.

\begin{remark}[Properties of the Evolution Equation]
\label{remark:charMethod}

The form of \eqref{eqn:def-elastic-evolution} suggests the method of characteristics.
Introducing the parametric variable $l$, we can write the solution of this PDE system through the solution of the following ODEs:
\begin{align}
    \label{eqn:our-char-t}
    &
    \deriv{t}{l} = 1, \quad t(l = 0) = 0 \qquad \Rightarrow \qquad t=l
    \\
    \label{eqn:our-char-x}
    &
    \deriv{x_i}{l} = v_i, \quad x_i(l=0) = c_i, \qquad i=1,2,3
    \\
    \label{eqn:our-char-F}
    &
    \deriv{(F_e)_{ij}}{l} = \parderiv{v_i}{x_k} (F_e)_{kj}, \quad (F_e)_{ij}\left(x_k(=c_k), l=t(=0)\right) = \text{initial or boundary condition} \qquad i,j,k=1,2,3
\end{align}
in a fixed Cartesian basis, and we have used index notation with implied summation.

The equations of the characteristic curves are defined by \eqref{eqn:our-char-t} and \eqref{eqn:our-char-x}:
\begin{equation}
    \hat{\bfx} = \left(\hat{x}_1(l,c_1,c_2,c_3) , \hat{x}_2(l,c_1,c_2,c_3) , \hat{x}_3 (l,c_1,c_2,c_3) , \hat{t}(l)\right)
\end{equation}
and \eqref{eqn:our-char-F} governs the behavior of $\bfF_e$ along the characteristic curves.
\end{remark}
\begin{remark}[Pathlines of the motion]
Comparing the transport equation \eqref{eqn:def-elastic-evolution} with the transport equation of $\bfchi^{-1}$ provided in \cite{kamrin2009eulerian}, the characteristic lines of both equations are the same and the characteristic curves of the latter are the pathlines of the motion. Therefore, the characteristic curves of \eqref{eqn:def-elastic-evolution} are the pathlines of the motion. We will use this property together with \eqref{eqn:bc-balance-mass} in the following examples to determine the shape of the body.
\end{remark}


\subsection{Summary of the Proposed Approach}
\label{sec:summary}

Our approach requires the solution of the balances of mass and momentum, and the transport of the elastic deformation.
In summary, we solve the 3 coupled equations given by:
\begin{equation}
\label{eqn:summary}
    \deriv{\ }{t} \begin{Bmatrix}
     \rho \\ \bfv \\ \bfF_e 
    \end{Bmatrix}
    =
    \begin{Bmatrix}
    -\rho \divergence\bfv \\ \frac{1}{\rho}\divergence \hat\bfsigma(\bfF_e) \\ (\nabla\bfv) \bfF_e 
    \end{Bmatrix}
\end{equation}
where $\deriv{\ }{t}$ is the material time derivative.


\section{A hollow sphere with accretion at the inner and ablation at the outer surface}
\label{sec:inne_acc}

In this section, we consider a hollow sphere with accretion and ablation. 
The accretion occurs on the inner fixed boundary and the ablation occurs on the outer boundary simultaneously. 
This example is motivated by the assembly and disassembly of actin micro-filaments during polymerization and depolymerization on a spherical bead \cite{pantaloni2001mechanism,noireaux2000growing}, and builds on the model of this phenomenon proposed by \cite{tomassetti2016steady}.
While we find the general solution, in the special case that there is no initial body, our solution reduces to that studied in \cite{tomassetti2016steady}, where it was solved using a Lagrangian approach with a 4-d manifold as the reference configuration. 

\begin{figure}[!ht]
	\centering
	\includegraphics[scale=0.75]{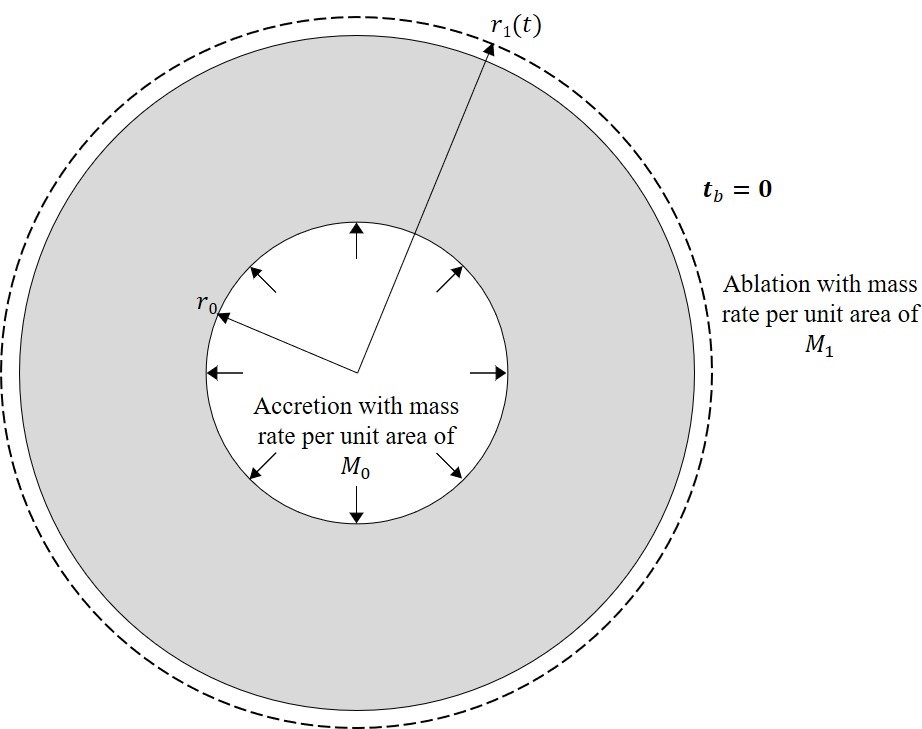}
	\caption{Schematic view of the hollow spherical body with accretion at inner fixed radius and ablation at the outer time dependent radius}
	\label{fig:inner_acc}
\end{figure}

For the stress response, we use an incompressible Neo-Hookean form:
\begin{equation}
    \label{eqn:const-law}
    \hat\bfsigma(\bfF) = -p \bfI + \parderiv{W(\bfF)}{\bfF} \bfF^T
\end{equation}
where $p$ is the pressure and $W$ is the stored energy density function.
    We note that the incompressible assumption violates the invertibility assumption of the stress-response function.
    For this particular problem, we do not use this assumption, and in the general setting, it is sufficient to have invertibility for the part of the stress derived from $W$.

The schematic view of the problem is depicted in Figure \ref{fig:inner_acc}.
We denote quantities related to the inner and outer surfaces using the subscript $0$ and $1$ respectively.
We work in the spherical coordinate system:
\begin{equation}
\label{eqn:coordinate_Sys_inner}
    \hat{\bfe}_r=\parderiv{\bfx}{r},\ \ \hat{\bfe}_\theta=\frac{1}{r}\parderiv{\bfx}{\theta}, \ \
    \hat{\bfe}_\phi=\frac{1}{r \sin\theta}\parderiv{\bfx}{\phi}
\end{equation}

\subsection{Boundary and Initial Conditions}

The inner radius is fixed and is equal to $r_0$, which means that the boundary velocity ${\bfV_b}_0$ in \eqref{eqn:bc-balance-mass} is equal to zero. 
Particles with given time-dependent velocity $\hat{\bfv}_{a_0} = - \dot{Z}_0(t) \hat{\bfe}_r, \ (\dot{Z}_0<0,\ Z_0(0)=0)$ approach the body at the inner surface and attach to it. 
Therefore, $M_0 = - \rho \dot{Z}_0(t) > 0$.
Then, using \eqref{eqn:bc-balance-mass}, it is clear that the normal velocity of the material particles right after attachment at $r=r_0$ is equal to $\dot{Z}_0(t)$, i.e. $\bfv = - \dot{Z}_0(t) \hat{\bfe}_r$.
This will provide a boundary condition for computing the velocity field later on.

Also, mass is removed from the outer surface with the ablation rate $M_1<0$. In the treadmilling regime, the boundary velocity is zero and the geometry of the body is constant, consequently $M_1 = - \rho \bfv (r=r_1(t)) \cdot \hat{\bfe}_r$.

The outer boundary is traction-free.

These velocity relations correspond to those in \cite{tomassetti2016steady}. 
However, they obtain $\frac{M_0}{\rho} = -\dot{Z}_0$ and $\frac{M_1}{\rho} = \dot{Z}_1(t) < 0$ using a sophisticated coupled analysis based on diffusion and kinetics; here, we will simply assume them as given.

We next turn to the deformation and stress states of the added particles.
Since the added particles are composed of the same material as the body, they have the same stress response.
    We assume that no tractions are imposed on the inner surface and the attaching particles are stress-free, up to an interior hydrostatic pressure $p_0$ that must be determined, i.e., $\bfsigma(r=r_0) = -p_0 \bfI$.
    From the assumption of incompressibility, we can choose $\bfF_e = \bfF = \bfFr = \bfI$.
The outer boundary is an outflow boundary and requires no boundary condition on $\bfF_e$.

Finally, for the initial conditions, we use $\bfF_e$ that is consistent with the boundary conditions and equilibrium equation.

\subsection{Reduced Governing Equations}

As the accretion velocity is radial and there is no shear traction, we assume that the motion is independent of $\theta$ and $\phi$.
The velocity $\bfv$ then simplifies to:
\begin{equation}
    \bfv = v_r(r,t) \hat{\bfe}_r
\end{equation}
Based on the definition of the problem, the governing equations \eqref{eqn:summary} together with the stress response \eqref{eqn:const-law} for a Neo-Hookean incompressible body reduce to the following :
\begin{align}
    &
    \label{eq:inner-acc-mass}
    \divergence(\bfv) = 
    \parderiv{v_r}{r} + \frac{2v_r}{r} = 0
    \quad \text{on }
    r_0 < r < r_1(t)
    \\
    &
    \label{eqn:bc-inner-vel-inner}
    \bfv \cdot \hat{\bfe}_r = \frac{M_0}{\rho} = - \dot{Z}_0 (t), 
	\quad
    \bfV_{b_0} = \mathbf{0}
    \quad \text{at } r=r_0
    \\
    &
    \label{eqn:bc-outer-vel-inner}
    \frac{M_1}{\rho} = \dot{Z}_1 (t) 
    \quad \text{at }
    r=r_1(t)
    \quad \Rightarrow \quad 
    \bfV_{b_1} \cdot \hat{\bfe}_r = v_r(r=r_1 , t) + \dot{Z}_1 (t)
    \\
    &
    \label{eqn:mom-inner}
    \divergence(\bfsigma) = \left( \parderiv{\sigma_{rr}}{r} + \frac{2}{r} (\sigma_{rr} - \sigma_{\theta \theta}) \right) \hat{\bfe}_r = \mathbf{0}
    \quad \text{on }
    r_0 < r < r_1(t) 
    \\
    &
    \label{eqn:bc-outer-mom-inner}
    \bft_b =  \mathbf{0}
    \quad \text{at }
    r=r_1(t)
    \\
    &
    \label{eqn:Fdot_inner}
    \parderiv{\ }{t}
    \begin{bmatrix}
	    F_{e_{rr}} & F_{e_{r \theta}} & F_{e_{r \phi}} \\
	    F_{e_{\theta r}} & F_{e_{\theta \theta}} & F_{e_{\theta \phi}} \\
	    F_{e_{\phi r}} & F_{e_{\phi \theta}} & F_{e_{\phi \phi}}
    \end{bmatrix}
    + v_r \parderiv{\ }{r}
    \begin{bmatrix}
	    F_{e_{rr}} & F_{e_{r \theta}} & F_{e_{r \phi}} \\
	    F_{e_{\theta r}} & F_{e_{\theta \theta}} & F_{e_{\theta \phi}} \\
	    F_{e_{\phi r}} & F_{e_{\phi \theta}} & F_{e_{\phi \phi}}
    \end{bmatrix}
    = 
    \begin{bmatrix}
	    \parderiv{v_r}{r} & 0 & 0 \\
	    0 & \frac{v_r}{r} & 0 \\
	    0 & 0 & \frac{v_r}{r}
    \end{bmatrix}
    \begin{bmatrix}
	    F_{e_{rr}} & F_{e_{r \theta}} & F_{e_{r \phi}} \\
	    F_{e_{\theta r}} & F_{e_{\theta \theta}} & F_{e_{\theta \phi}} \\
	    F_{e_{\phi r}} & F_{e_{\phi \theta}} & F_{e_{\phi \phi}}
    \end{bmatrix}
    \\
    &
    \label{eqn:Fdot_bc_ic_inner}
    \bfF_e (r , t=0) \text{ given}
    , \quad
    \bfF_e (r = r_0 , t) =\bfI 
    \\
    &
    \label{eqn:neo-Hookean-inner}
    \hat\bfsigma(\bfF_e) =
    -p \bfI + G \bfF_e \bfF_e^T, \quad
    G \text{ is the shear modulus}
\end{align}

\subsection{Solution}

Using the boundary condition \eqref{eqn:bc-inner-vel-inner}, the continuity equation \eqref{eq:inner-acc-mass} gives us the velocity field:
\begin{equation}
\label{eq:inner-acc-vel}
    \bfv = -\dot{Z}_0 \frac{ r_0^2}{r^2} \hat{\bfe}_r
\end{equation}
We next notice that we can write \eqref{eqn:Fdot_inner} as:
\begin{equation*}
    \left( \parderiv{\ }{t} + v_r \parderiv{\ }{r} - A \right) F_{e_{ij}}=0
    , \quad
    A = 
    \begin{cases}
                \parderiv{v_r}{r}
               & \text{if} \ \ ij = rr, r\theta, r\phi \ \  \\
               \\
                \frac{v_r}{r} 
                & \text{otherwise}  \ \ 
    \end{cases}
\end{equation*}
These differential equations are homogeneous, i.e. no source term on the right, and therefore those components with zero initial and boundary conditions remain zero during the growth process. 
Hence, the off-diagonal components of $\bfF_e$ are zero throughout the growth process, assuming that $\bfF_e$ is diagonal at the initial time.
We can then represent $\bfF_e$ as:
\begin{equation}
\label{eqn:def-grad-ex1}
    \bfF_e = 
    F_{e_{rr}} \hat{\bfe}_r \otimes \hat{\bfe}_r 
    + F_{e_{\theta\theta}} \hat{\bfe}_\theta \otimes \hat{\bfe}_\theta + 
    F_{e_{\phi\phi}} \hat{\bfe}_\phi \otimes \hat{\bfe}_\phi
\end{equation}

The evolution equation \eqref{eqn:Fdot_inner} and the initial (assuming $F_{e_{\theta \theta}} (t=0) = F_{e_{\phi \phi}} (t=0)$) and boundary conditions \eqref{eqn:Fdot_bc_ic_inner} are the same for $F_{e_{\theta\theta}}$ and $F_{e_{\phi\phi}}$, as we would expect from spherical symmetry, so we solve only for $F_{e_{\theta\theta}}$. 
Thus, only two of the nine scalar equations for the components of $\bfF_e$ are independent and non-trivial.
Finally, after substituting the velocity \eqref{eq:inner-acc-vel}, we obtain two decoupled equations:
\begin{align}
\label{eqn:Frr-inner}
    & \parderiv{F_{e_{rr}}}{t} - \dot{Z}_0 \frac{r_0^2}{r^2} \parderiv{F_{e_{rr}}}{r} = 2 \dot{Z}_0 \frac{r_0^2}{r^3} F_{e_{rr}}
    \\
    \label{eqn:Ftt-inner}
    & \parderiv{F_{e_{\theta\theta}}}{t} - \dot{Z}_0 \frac{r_0^2}{r^2} \parderiv{F_{e_{\theta\theta}}}{r} = -\dot{Z}_0 \frac{r_0^2}{r^3} F_{e_{\theta\theta}}
\end{align}

\subsubsection{The Characteristic Curves}

To solve \eqref{eqn:Frr-inner} and \eqref{eqn:Ftt-inner}, we use the method of characteristics, similar to the discussion in remark \ref{remark:charMethod}, because of the convective nature of the equations. 
The parametric representation of the characteristic lines for both of these equations is:
\begin{align}
    &
    \label{eqn:char-t-ex1}
    \deriv{t}{l} = 1 
    \ \ \ \Rightarrow \ \ \ 
    \hat{t} (l) = l + C_0
    \\
    &
    \label{eqn:char-r-ex1}
    \deriv{r}{l} = -\dot{Z}_0 \frac{r_0^2}{r^2}
    \ \ \ \Rightarrow \ \ \ 
    \hat{r}^3(l) = -3 Z_0(\hat{t}) r_0^2 + C_1
\end{align}
where $l$ is the parametric variable along the characteristic curves.

\begin{figure}[htb!]
    \centering
    \includegraphics[scale=0.25]{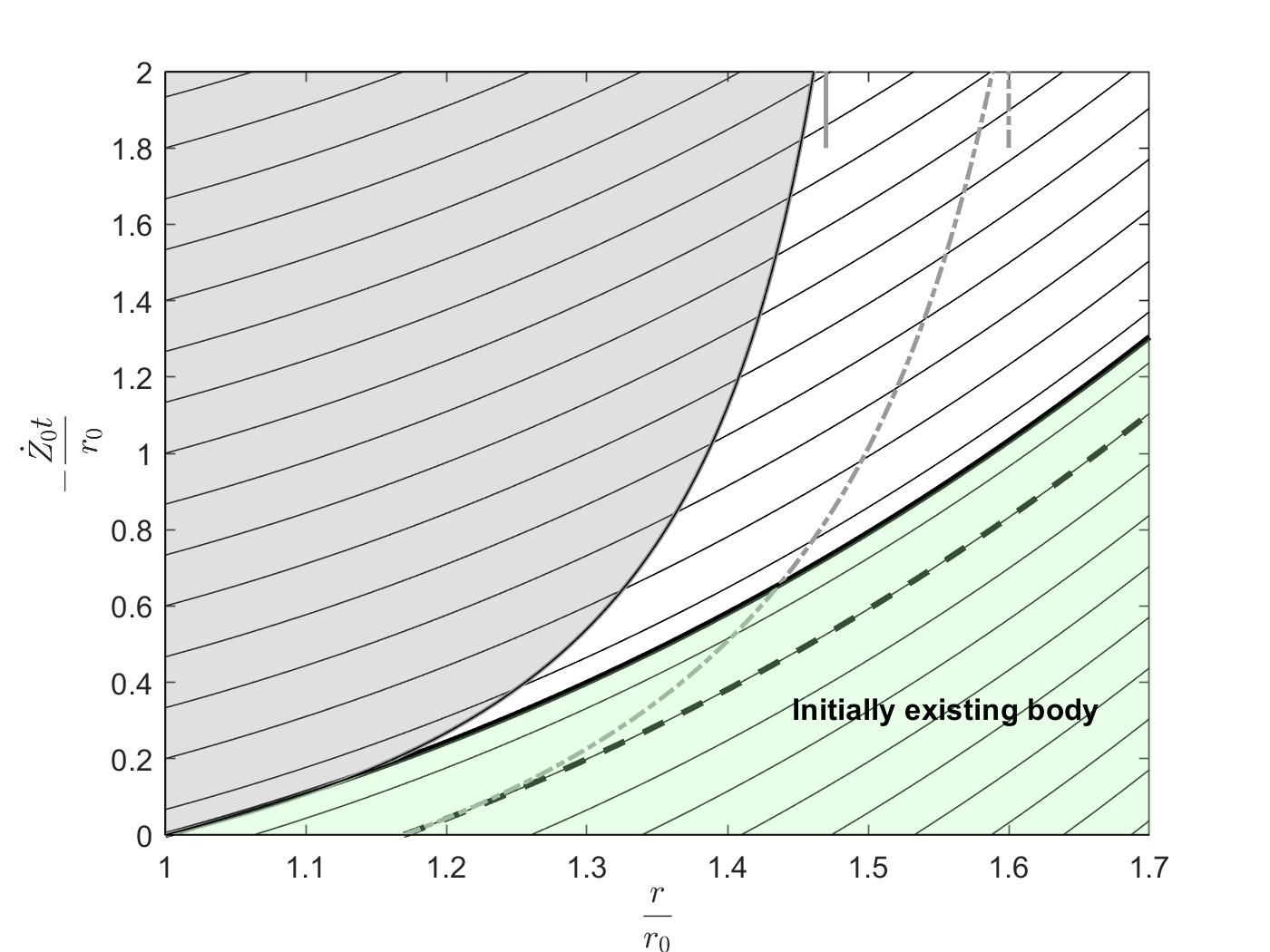}
    \caption{Characteristic curves which convect initial and boundary deformation through out the space-time domain. 
    (1) The solid non-bold curves are the characteristic curves, 
    (2) the bold black curve distinguishes the initial body (green region) from the added particles (non-green region), and is also the outer radius when there is no ablation and no initial body, 
    (3) the black dashed curve shows the outer radius in the case of no ablation and existence of an initial body, 
    (4) the gray shaded area shows the body when there is ablation and no initial body, and
    (5) the gray dot-dash curve shows the outer radius in the case of both ablation and existence of the initial body.}
    \label{fig:inner_charac}
\end{figure}

Assuming that all the characteristic curves start from $t=0$, we have for all curves that $C_0 = 0$ and $\hat{t}=l$, although the information about $\bfF_e$ might be provided at the middle of the curve (at the attachment time $\tau$ of each particle).
However, different characteristic curves have different initial locations $\hat{r}(l=0)$, so $C_1$ is the quantity that differentiates characteristic curves. 
For the special case of constant $\dot{Z}_0$, different characteristic curves are plotted in Figure \ref{fig:inner_charac}. The information about $\bfF_e$ is carried along these characteristic curves which are also pathlines of the motion.

Instead of distinguishing different characteristic curves by their intersection with the line $t=0$, one could alternatively determine them by their attachment time $\tau$ to the body. 
In other words, assuming that all characteristic curves start from $r=r_0$ at $l=0$, then $C_1 = r_0^3$ for all characteristic curves but they are distinguished by parameter $C_0$ which is the attachment time $\tau$ that was defined by \cite{skalak1997kinematics} and \cite{tomassetti2016steady}. 
We note that the extrapolation of characteristic curves to the $t=0$ axis and representing them by the attachment time $\tau$ provides a 4-d manifold that is denoted the ``placement map'' in \cite{tomassetti2016steady}.
In the case that there is no initial body at $t=0$ and all the particles were added to the body via accretion, these two choices of distinguishing characteristic curves are equivalent from the continuum mechanics perspective.
However, if there is a stressed or an unstressed initial body at $t=0$ with accretion and ablation occurring at its boundaries, the approach of using the attachment time cannot as easily predict the deformation of the body during the growth process.
In other words, there is no distinction between the particles that were initially in the body (i.e., with $\tau=0$) at different radial locations, so their motion can only be considered as a rigid body motion at any later time as in \cite{skalak1997kinematics} or it needs extra effort.

The information about $\bfF_e$ is carried along these characteristic curves. 
It can be seen from figure \ref{fig:inner_charac} that, based on the origin of each curve, they may carry either initial condition or boundary condition information. 
So the space-time domain splits into two regions: (1) the region influenced by the initial condition (green region), and (2) the region influenced by the boundary condition (non-green region). 
The deformation gradient is affected by the initial condition in the lower region of the black bold curve (information comes from the line $t=0$), and by the boundary condition in the upper region of the black bold curve (information comes from the line $\frac{r}{r_0}=1$).

Each characteristic curve in the lower region shows the location of a particle which was at specific non-dimensional radius $\frac{r}{r_0}$ at $t=0$. 
Similarly, each characteristic curve in the upper region shows the location of a particle which was added to the body at time $\tau$; the non-dimensional time $- \frac{\dot{Z}_0 \tau}{r_0}$ is given by the intersection of the characteristic curve with the line $\frac{r}{r_0}=1$. 
It should be noted that these pathlines show the motion of the particles inside the growing body, so they are only valid within the growing body though we have extended them outside the body.
The time-dependent geometry of the growing body depends on the initial shape of the body and the boundary velocity; the latter depends on the rate of accretion or ablation and the particle velocity at the boundary, or it might be given as in the case of the inner fixed boundary. 

When there is no ablation ($M_1=0$) and the growth starts from a substrate with radius $r_0$, i.e. there is no body at $t=0$, the black solid bold curve in figure \ref{fig:inner_charac} shows the outer radius of the body at each time because only the motion of the particles at the outer boundary determines the velocity and location of the boundary.
Then, the time-dependent growing body is a hollow sphere with a constant non-dimensional inner radius of one and non-dimensional outer radius given by the black solid bold curve. 
The distinguishing constant $C_1$ in \eqref{eqn:char-r-ex1} for the outer radius in the special case of no ablation and no initial body is obtained by substituting $\hat{r}=r_0, \hat{t} = l = 0$ and solving for $C_1$, which gives
\begin{equation}
    \label{eqn:inner-front}
    r^3 = r_0^3 - 3 r_0^2 Z_0(t)
\end{equation}
This is the equation of the black solid bold curve.

However, if there exists an initial body with inner and outer radius $r_0$ and $r_1(t=0)$, and there is no ablation, the outer radius of the hollow sphere is defined by the characteristic line that was at $r=r_1(t=0)$ at $t=0$ (dashed bold black curve). 
For the general case with ablation and an initial body, the outer radius of the growing body (denoted by $r_1(t)$ in this problem) increases, but at a smaller rate than the case with no ablation, and the rate depends directly on the given function $\dot{Z}_1$. 
In this case, the outer radius is depicted (schematically) by the dot-dash gray curve in figure \ref{fig:inner_charac}.
If this radius reaches a steady-state constant value, the body will be ``treadmilling'' \cite{tomassetti2016steady}.

Finally, we consider the case when there is ablation but no initial body, following \cite{tomassetti2016steady}.
The non-dimensional outer radius of $\frac{r_1(t)}{r_0}$ which depends on $\dot{Z}_1$ is shown (schematically) by the solid gray curve and the growing body is shown by the gray area in figure \ref{fig:inner_charac}. 
As the inner radius is fixed, the lower limit of the growing body is the vertical line at non-dimensional radius $\frac{r}{r_0} = 1$, and the upper limit of the area is the non-dimensional outer radius $\frac{r_1(t)}{r_0}$ at each time (the gray curve). 
As can be seen from the figure, $\bfF_e$ in this case is only affected by the boundary conditions, since there can be no initial conditions absent an initial body.

\subsubsection{Evolution of the deformation gradient}

To compute the nonzero components of deformation gradient, we consider the scalar equations \eqref{eqn:Frr-inner} and \eqref{eqn:Ftt-inner} separately as their right hand sides are different.

\begin{proof}[Solving \eqref{eqn:Frr-inner} for $F_{e_{rr}}$.]
    The differential equation of $F_{e_{rr}}$ along characteristic curves which is parametrized with $l$ is:
    \begin{equation*}
        \deriv{F_{e_{rr}}}{l} = 2\dot{Z}_0 \frac{r_0^2}{r^3}F_{e_{rr}}
    \end{equation*}
    Combining this with \eqref{eqn:char-r-ex1} gives:
    \begin{equation*}
        \frac{\dm F_{e_{rr}}}{F_{e_{rr}}} = -2 \frac{\dm r}{r}
        \quad \Rightarrow\quad
        F_{e_{rr}}=\frac{f(C_1)}{r^2}=\frac{f \left(r^3+3Z_0r_0^2 \right)}{r^2}
    \end{equation*}
    The unknown function $f$ is the constant of integration which is different for each characteristic curve, and it is a function of $C_1$ which distinguishes different characteristic curves.
    It is to be determined using initial and boundary conditions.
    Consider the initial condition of $F_{e_{rr}} (r,t=0) = F_{e_{rr_0}}(r)$, so $f$ is as follows for the region influenced by the initial condition:
    \begin{equation*}
        F_{e_{rr_0}}(r) = \frac{f(r^3)}{r^2} 
        \quad\Rightarrow\quad
        f(x) = x^{2/3} F_{e_{rr_0}}(x^{1/3}) 
    \end{equation*}
    The boundary condition is also $F_{e_{rr}} (r = r_0,t) = 1$, so $f$ is as follows for the region influenced by the boundary condition:
    \begin{equation*}
        1 = \frac{f \left(r_0^3+3Z_0(t)r_0^2\right)}{r_0^2} 
        \quad\Rightarrow\quad
        f(x) = \const = r_0^2
    \end{equation*}
    The function $f$ is constant because the argument depends on $t$ but it is equal to constant value $r_0^2$.
    
    Putting these together, we have:
    \begin{equation*}
        F_{e_{rr}}(r,t)= 
            \begin{cases}
               \frac{(r^3+3 Z_0 r_0^2) ^ \frac{2}{3} F_{e_{rr_0}} \Big((r^3+3 Z_0 r_0^2) ^ \frac{1}{3} \Big)}{r^2}
               & \text{if} \ \ r^3 < r_0^3 - 3 r_0^2 Z_0(t) 
               \\
                \left(\frac{r_0}{r}\right)^2 
                & \text{if} \ \ r^3 > r_0^3 - 3 r_0^2 Z_0(t)
            \end{cases}
    \end{equation*}
    Therefore, in the problem defined in \cite{tomassetti2016steady} where there is no body at initial state but there is ablation at the outer surface, $F_{e_{rr}}$ is $\left(\frac{r_0}{r}\right)^2$ at all times.
\end{proof}

\begin{proof}[Solving \eqref{eqn:Ftt-inner} for $F_{e_{\theta \theta}}$.]
    The differential equation for $F_{e_{\theta \theta}}$ along the characteristic curves which are parametrized by $l$ is:
    \begin{equation*}
        \deriv{F_{e_{\theta \theta}}}{l} = - \dot{Z}_0 \frac{r_0^2}{r^3}F_{e_{\theta \theta}}
    \end{equation*}
    Combining this with \eqref{eqn:char-r-ex1} gives:
    \begin{equation*}
        \frac{\dm F_{e_{\theta \theta}}}{F_{e_{\theta \theta}}} = \frac{\dm r}{r}
        \quad\Rightarrow\quad
        F_{e_{\theta \theta}}=r g(C_1)=r g \left(r^3+3Z_0r_0^2 \right)
    \end{equation*}
    The unknown function $g$ is the constant of integration which is different for each characteristic curve, so it is a function of $C_1$ which distinguishes different characteristic curves. 
    Using the initial condition of $F_{e_{\theta \theta}} (r,t=0) = F_{e_{\theta \theta_0}}(r)$ and boundary condition, similar to the approach used in computing $F_{e_{rr}}$, we have:
    \begin{equation*}
        F_{e_{\theta \theta}}(r,t)= 
            \begin{cases}
               \frac{r F_{e_{\theta \theta_0}}\Big( (r^3+3 Z_0 r_0^2) ^ \frac{1}{3} \Big)}{(r^3+3 Z_0 r_0^2) ^ \frac{1}{3}}
               & \text{if} \ \ r^3 < r_0^3 - 3 r_0^2 Z_0(t)
               \\
                \frac{r}{r_0} 
                & \text{if} \ \ r^3 > r_0^3 - 3 r_0^2 Z_0(t)
            \end{cases}
    \end{equation*}
    For the problem defined in \cite{tomassetti2016steady}, $F_{e_{\theta \theta}}$ is $\frac{r}{r_0}$ at all times.
\end{proof}

\subsection{Solution for the special case with no initial body and with ablation at the outer surface}

For the example defined in \cite{tomassetti2016steady} where there is no initial body and $r^3 > r_0^3 - 3 r_0^2 Z_0(t)$, our calculation of $\bfF_e$ is:
\begin{equation}
\label{eqn:finalF-inner}
    \bfF_e = 
    \left(\frac{r_0}{r}\right)^2 \hat{\bfe}_r \otimes \hat{\bfe}_r + \frac{r}{r_0} \left(\hat{\bfe}_\theta \otimes \hat{\bfe}_\theta + \hat{\bfe}_\phi \otimes \hat{\bfe}_\phi\right)
\end{equation}
As the process of addition of material is smooth, there is no jump in traction boundary conditions to change $\bfF_e$ instantaneously, and the added particles and the traction boundary condition are in equilibrium.
Therefore, $\bfF_{relax} = \bfI$.

Our calculation of $\bfF$ matches with that computed in \cite{tomassetti2016steady}; a superficial difference is that they represent the deformation gradient in the mixed basis with the first leg in the current and the second leg in the reference.
Further, our calculation of the equations of motion (characteristic curves) derived in \eqref{eqn:char-r-ex1} matches with the current location of particles that was added to the body at time $t_0(Z)$ derived in \cite{tomassetti2016steady} with $C_1 = r_0^3 + 3r_0^2 Z$, where $t_0(Z)$ is the attachment time of a particle that was at location $Z$ in the reference configuration at the time of attachment. 
Therefore, essentially, both $Z$ and $\hat{r}(l=0)$ are the same and the unit vectors in $r$ and $Z$ direction are identical. 
Further, the problem is independent of $\theta$ and $\phi$ and the coordinates are orthonormal, hence the coordinate vectors defined here and \cite{tomassetti2016steady} are the same, implying that the representations of $\bfF$ are also identical.

Finally, substituting \eqref{eqn:finalF-inner} in \eqref{eqn:neo-Hookean-inner}, the Cauchy stress tensor is:
\begin{equation*}
    \bfsigma = 
    \left(G \left(\frac{r_0}{r}\right)^4-p \right) \hat{\bfe}_r \otimes \hat{\bfe}_r + 
    \left(G \left(\frac{r}{r_0}\right)^2-p \right) \left(\hat{\bfe}_\theta \otimes \hat{\bfe}_\theta + \hat{\bfe}_\phi \otimes \hat{\bfe}_\phi\right)
\end{equation*}
To determine the unknown Lagrange multiplier function $p$, the balance of linear momentum \eqref{eqn:mom-inner} and the boundary condition \eqref{eqn:bc-outer-mom-inner} are used. Thus, the stress field in the growing body is:
\begin{align}
\label{eqn:finalSigma-inner}
    \frac{\bfsigma}{G} = 
        & \left[ \frac{1}{2} \left(\left(\frac{r_0}{r}\right)^4 -   \left(\frac{r_0}{r_1(t)}\right)^4 \right) + \left(\frac{r}{r_0}\right)^2 - \left(\frac{r_1(t)}{r_0}\right)^2 \right] \hat{\bfe}_r \otimes \hat{\bfe}_r + \nonumber
        \\
        & \left[ -\frac{1}{2} \left(\left(\frac{r_0}{r}\right)^4 + \left(\frac{r_0}{r_1(t)}\right)^4 \right) + 2 \left(\frac{r}{r_0}\right)^2 - \left(\frac{r_1(t)}{r_0}\right)^2 \right] \left(\hat{\bfe}_\theta \otimes \hat{\bfe}_\theta + \hat{\bfe}_\phi \otimes \hat{\bfe}_\phi\right)
\end{align}
which agrees with the Cauchy stress tensor derived in \cite{tomassetti2016steady}.


\section{A hollow cylinder with accretion at the outer surface}
\label{sec:outer_acc}

In the previous example, accretion occurred at a surface with a known location.
Here we will consider the case that accretion occurs at an unknown time-dependent surface.
Specifically, we consider a long hollow cylinder made of a general incompressible material with time-dependent inner and outer radii, and accretion occurring on the outer surface and traction applied on the inner surface.
The formulation is motivated by the wound roll process of additive manufacturing which has been applied to processes ranging from the packaging of flexible films in paper and magnetic tape industries \cite{yagoda1980resolution} to tissue engineered blood vessels \cite{konig2009mechanical}.
The stress response corresponds to an incompressible neo-Hookean material as in \eqref{eqn:const-law}.

We assume that we can neglect the end effects and assume plain strain.
Therefore, we use a 2-d polar coordinate system:
\begin{equation}
    \label{eqn:coorsys_polar_nd}
    \hat{\bfe}_r=\parderiv{\bfx}{r},\ \ 
    \hat{\bfe}_\theta=\frac{1}{r} \parderiv{\bfx}{\theta}, \ \
\end{equation}

This problem was solved in \cite{sozio2017nonlinear} using a geometrical Lagrangian approach with two different choices of metric for the initial and added parts of the body.
We discuss below that the choice of the growth velocity does not change the Eulerian formulation, but it can make the Lagrangian formulation more complicated as the time of attachment has to be related to the location in the reference configuration. 
We also show that one can solve an Eulerian growth process using the evolution of inverse of the motion via the method discussed in \cite{kamrin2009eulerian}, rather than computing the deformation gradient directly using \eqref{eqn:def-elastic-evolution}.

\subsection{Boundary and Initial Conditions}

We denote the time-dependent inner and outer radii by $r_{in}(t)$ and $r_{out}(t)$, their initial values by $r_{in}(t=0) = R_0$ and $r_{out}(t=0)=R_1$, and the radial velocities by $\dot{r}_{out}$ and $\dot{r}_{in}$, respectively. 

\begin{figure}[htb!]
    \centering
    \includegraphics[scale=0.8]{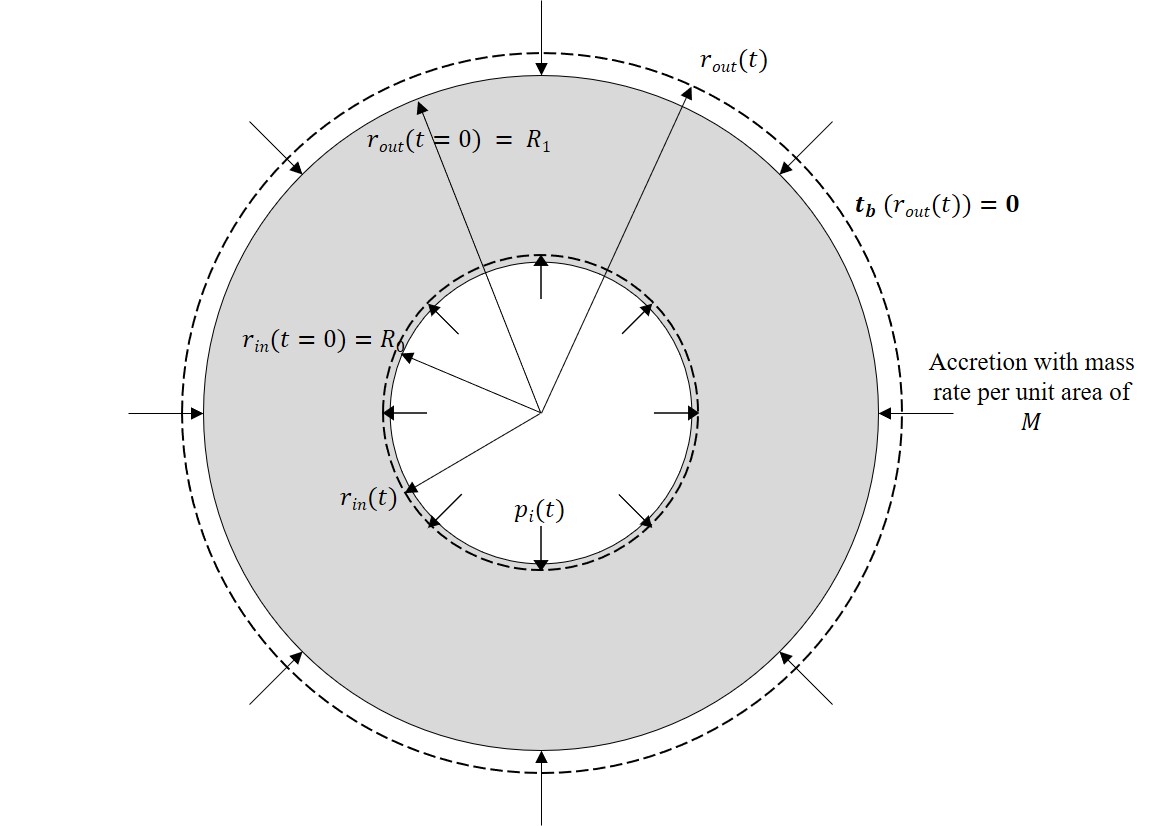}
    \caption{Schematic view of the  very long hollow cylinder body with outer accretion, the gray body is initial body and the dashed lines show the time-dependent boundary of the growing body}
    \label{fig:outer_acc_prob}
\end{figure}

The body is initially stress-free, and an internal pressure of $p_i(t)$ ($p_i(t=0) = 0$) applies a traction at the inner time-dependent surface of the cylinder, shown in figure \ref{fig:outer_acc_prob}. 

There is no ablation, and the accretion occurs only at the outer surface. 
The mass rate of growth $M$ at the outer surface per unit area is equal to $\rho  u_g (t)$. The added particles are stress-free and the growing surface is traction-free.
We aim to find the time-dependent geometry of the body, i.e. the expression for $r_{in}(t)$ and $r_{out}(t)$, and the stress distribution in the body.

\subsection{Reduced Governing Equations}

From the symmetry of the problem, we assume that no quantities of interest depend on $\theta$, and that the velocity $\bfv$ is radial and can be written $\bfv=v_r(r,t) \hat{\bfe}_r$.
Therefore, the governing equations \eqref{eqn:summary} together with \eqref{eqn:const-law}, written in the coordinate system defined in \eqref{eqn:coorsys_polar_nd}, reduces to the following:
\begin{align}
    &
    \label{eq:outer-acc-mass}
    \divergence(\bfv)=
    \parderiv{v_r}{r}+\frac{v_r}{r}
    = 0
    \quad \text{on} \quad
    r_{in}(t) < r < r_{out}(t)
    \\
    &
    \label{eqn:bc-outer-vel-outer}
    \dot{r}_{out} = \bfV_b \cdot \hat{\bfe}_r = 
    v_r (r=r_{out}(t) , t) + u_g(t)
    \quad \text{at} \quad
    r=r_{out}(t)
    \\
    &
    \label{eqn:bc-inner-vel-outer}
    \dot{r}_{in} = \bfV_b \cdot \hat{\bfe}_r = v_r (r=r_{in}(t) , t)
    \quad \text{at} \quad
    r=r_{in}(t)
    \\
    &
    \label{eqn:mom-outer}
    \divergence(\bfsigma) = \left(  \parderiv{\sigma_{rr}}{r}+\frac{1}{r}(\sigma_{rr}-\sigma_{\theta\theta}) \right) \hat{\bfe}_r = \mathbf{0}
    \quad \text{on} \quad
    r_{in}(t) < r < r_{out}(t) 
    \\
    &
    \label{eqn:bc-outer-mom-outer}
    \bft_b =  \mathbf{0}
    \quad \text{at} \quad
    r=r_{out}(t)
    \\
    &
    \label{eqn:bc-inner-mom-outer}
    \bft_b =  p_i(t) \hat{\bfe}_r
    \quad \text{at} \quad
    r=r_{in}(t)
    \\
    &
    \label{eqn:Fdot_outer}
    \parderiv{\ }{t}
    \begin{bmatrix}
    F_{e_{rr}} & F_{e_{r \theta}} \\
    F_{e_{\theta r}} & F_{e_{\theta \theta}} \\
    \end{bmatrix}
    + v_r \parderiv{\ }{r}
    \begin{bmatrix}
    F_{e_{rr}} & F_{e_{r \theta}} \\
    F_{e_{\theta r}} & F_{e_{\theta \theta}} \\
    \end{bmatrix}
    = 
    \begin{bmatrix}
    \parderiv{v_r}{r} & 0 \\
    0 & \frac{v_r}{r} \\
    \end{bmatrix}
    \begin{bmatrix}
    F_{e_{rr}} & F_{e_{r \theta}} \\
    F_{e_{\theta r}} & F_{e_{\theta \theta}} \\
    \end{bmatrix}
    \\
    &
    \label{eqn:Fdot_bc_ic_outer}
    \bfF_e = \bfI \text{ at } t=0, \ R_0<r<R_1, 
    \quad \text{and} \quad
    \bfF_e=\bfI \text{ at } r=r_{out}(t)
    \\
    \begin{split}
    \label{eqn:const-law-outer}
        &
        \hat\bfsigma(\bfF_e) =
        -p \bfI + \parderiv{W(\bfF_e)}{\bfF_e} \bfF_e^T =
        \left(-p + 2 I_2 \parderiv{W}{I_2} \right)\bfI + 2\parderiv{W}{I_1} \bfF_e \bfF_e^T - 2 \parderiv{W}{I_2} (\bfF_e \bfF_e^T)^{-1}
        , \quad 
        \\
        &
        \qquad\qquad
        \text{where } I_1 = \trace(\bfB_e) \text{ and } I_2 = \frac{\trace (\bfB_e)^2 - \trace (\bfB_e^2) }{2} \text{ are invariants of } \bfB_e = \bfF_e \bfF_e^T
    \end{split}
\end{align}
As the initial and inlet boundary conditions of $\bfF_e$ are consistent, there is no jump in traction boundary conditions to change $\bfF_e$ instantly.
Also, the added particles and the traction boundary conditions have consistent stresses. 
So, we also keep the boundary and initial conditions for $\bfF$ equal to $\bfI$ and the relaxed shape of the added particles is equal to their shape at the time of attachment, then $\bfF_{relax} = \bfI$.

\subsection{Solution}

\begin{figure}[htb!]
    \centering
    \includegraphics[scale=0.25]{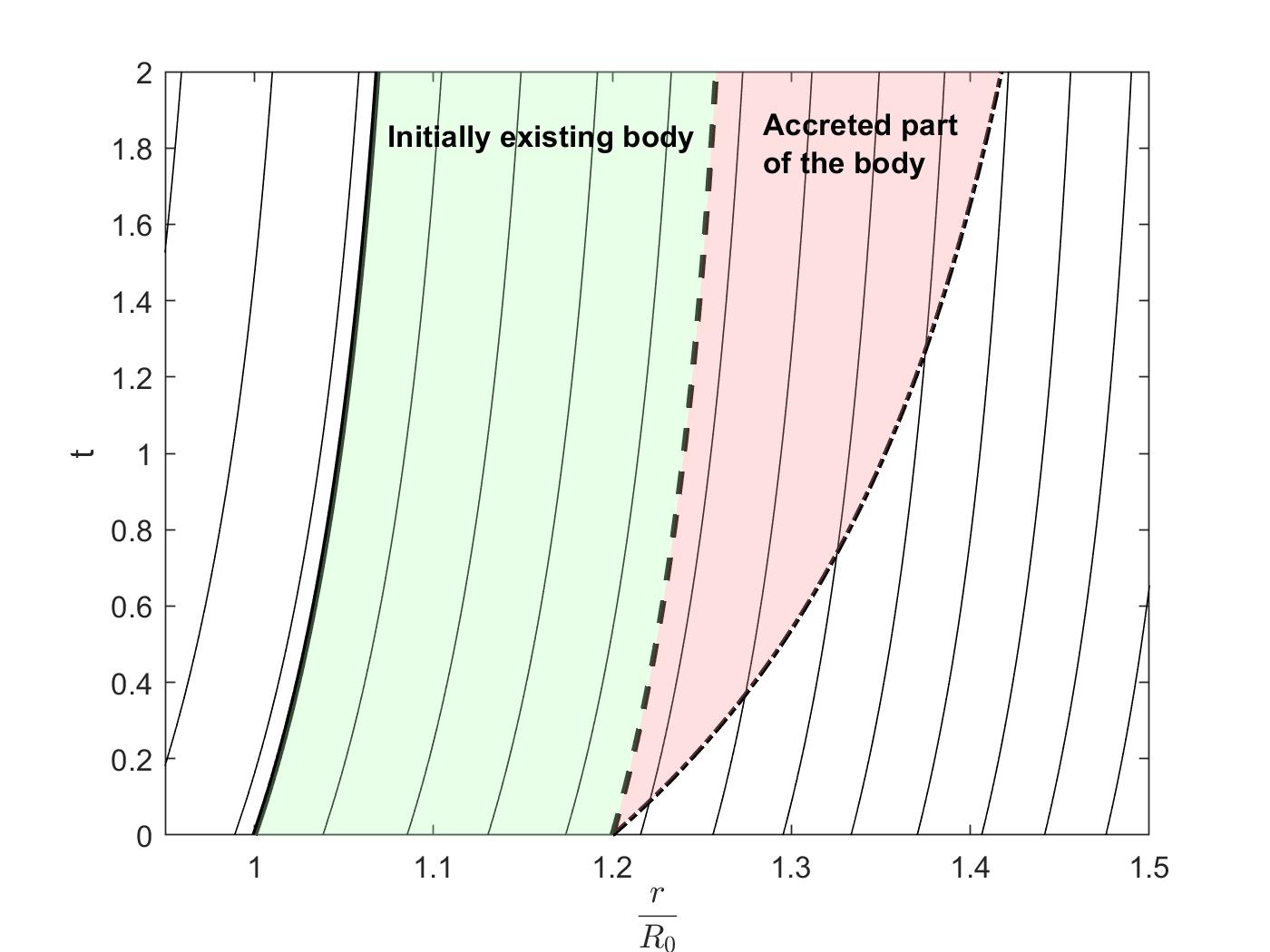}
    \caption{
        Schematic view of pathlines of the motion in space-time.
        The outer radius is marked by the bold dot-dash line, and the inner radius by the solid bold line; the time-dependent growing body lies between these radii.
        The characteristic curve starting at $r=R_0$ at $t=0$ is shown by the solid bold line.
        The colors mark the space-time regions corresponding to the initial body and the accreted mass.
    }
    \label{fig:outer_acc_charac}
\end{figure}

The solution of the continuity equation \eqref{eq:outer-acc-mass} using the boundary condition of \eqref{eqn:bc-inner-vel-outer} is:
\begin{equation}
    \label{eqn:vel_outer}
    \bfv=\frac{r_{in}}{r}  \dot{r}_{in} \hat{\bfe}_r
\end{equation}
and $r_{in}(t)$ and consequently $\dot{r}_{in}(t)$ will be computed later.
Substituting this in \eqref{eqn:bc-outer-vel-outer}, we obtain a differential equation that can be later solved for the outer radius:
\begin{equation}
    \label{eqn:in_out_vel_relation_outeracc}
    r_{in} \dot{r}_{in} = r_{out} \dot{r}_{out} - r_{out} u_g(t), \quad
    \text{ with } r_{out}(t=0)=R_1 \text{ and } r_{in}(t=0)=R_0
\end{equation}
Similar to the previous example, $\nabla \bfv$ is diagonal, so the four PDEs for the evolution of the components of $\bfF_e$ are decoupled and homogeneous.  
As both the boundary and initial conditions for $\bfF_e$ are identity, the off-diagonal components of $\bfF_e$ remain zero during the growth process.
Using \eqref{eqn:Fdot_outer} and \eqref{eqn:vel_outer}, the evolution equations for the diagonal components of $\bfF_e$ are:
\begin{align}
\label{eqn:Frr-outer}
    & \parderiv{F_{e_{rr}}}{t} + \dot{r}_{in} \frac{r_{in}}{r} \parderiv{F_{e_{rr}}}{r} = - \dot{r}_{in} \frac{r_{in}}{r^2} F_{e_{rr}}
    \\
    \label{eqn:Ftt-outer}
    & \parderiv{F_{e_{\theta\theta}}}{t} + \dot{r}_{in} \frac{r_{in}}{r} \parderiv{F_{e_{\theta\theta}}}{r} = \dot{r}_{in} \frac{r_{in}}{r^2} F_{e_{\theta\theta}}
\end{align}
We solve these using the method of characteristics, similar to the previous example.
The parametric equation of the characteristic curves are:
\begin{align}
    \label{eqn:char-line-t-ex2}
   & \deriv{t}{l} = 1 
   \Rightarrow
   \hat{t} = l + C_0
    \\
    \label{eqn:char-line-r-ex2}
   &\deriv{r}{l} = {\dot{r}_{in}\frac{r_{in}}{r}} 
\end{align}
As in the previous example, we assume that all the characteristic curves start at $t=0$, so $\hat{t}(l=0)=0 \Rightarrow C_0 = 0$ for all characteristic curves. 
Further, substituting $\dm l = \dm t$ in \eqref{eqn:char-line-r-ex2} and using the identity $\dot{r}_{in} r_{in} = \frac{1}{2} \deriv{\ }{t}r_{in}^2$, the function $\hat{r}(l)$ along the characteristic curves is:
\begin{equation}
\label{eqn:outeracc_char_in}
   r \dm r = \frac{1}{2} \dot{(r_{in}^2)} \dm t
   \quad
   \Rightarrow
   \quad
   \hat{r}^2(l) + C_1 =
   r_{in}^2(\hat{t}(l))
\end{equation}
The integration constant $C_1$ distinguishes different characteristic curves. 
Obviously, the pathline with $C_1=0$ shows the location of the inner boundary. 
However, the equation of $r_{in}(t)$ has not been computed yet, and will be done further below using the traction boundary conditions.
Further, given $r_{in}(t)$, we can compute $r_{out}(t)$ using \eqref{eqn:in_out_vel_relation_outeracc}.
Figure \ref{fig:outer_acc_charac} shows a schematic view of the pathlines and body shape.

We can also substitute for $\dot{r}_{in}r_{in}$ in \eqref{eqn:outeracc_char_in} from \eqref{eqn:in_out_vel_relation_outeracc} to find the equation for the characteristic curves:
\begin{equation}
    \label{eqn:outeracc_char_out}
    \hat{r}^2(l) + C_2 
    = 
    r_{out}^2(\hat{t}(l)) - 2 \int_{0}^{\hat{t}(l)} u_g(t) r_{out}(t) \dm t
\end{equation}
where $C_2$ is a constant of integration, and is different for different characteristic curves.

Using \eqref{eqn:Frr-outer} and \eqref{eqn:Ftt-outer}, $F_{e_{rr}}$ and $F_{e_{\theta\theta}}$ along the characteristic curves are governed by:
\begin{equation*}
   \deriv{F_{e_{rr}}}{l}={-\dot{r}_{in}\frac{r_{in}}{r^2}F_{e_{rr}}},
   \quad 
   \deriv{F_{e_{\theta \theta}}}{l}={\dot{r}_{in}\frac{r_{in}}{r^2}F_{e_{\theta \theta}}} 
\end{equation*}
Combining this with \eqref{eqn:outeracc_char_in}, the general form of $F_{rr}$ and $F_{\theta\theta}$ are:
    \begin{align}
        \label{eqn:Frr_general_outer}
       & \hat{r}(l,C) F_{e_{rr}} = g_1(C)
       \\
       \label{eqn:Ftt_general_outer}
       & \frac{F_{e_{\theta\theta}}}{\hat{r}(l,C)}=g_2(C)
    \end{align}
where $g_1$ and $g_2$ are constants of integration which are different for each characteristic line.
Therefore $g_1$ and $g_2$ are functions of either $C_1$ or $C_2$, depending on whether we use \eqref{eqn:outeracc_char_in} or \eqref{eqn:outeracc_char_out} to distinguish the different characteristic curves. 

Considering the characteristic curve starting at $r=R_1$ at $t=0$ (marked by a dashed bold line in figure \ref{fig:outer_acc_charac}), we substitute $(r=R_1, l=0)$ to find $C_1$ in \eqref{eqn:outeracc_char_in}, and hence the equation of this characteristic curve is:
\begin{equation}
    \label{eqn:outer_rin_r(R_0)}
    R_0^2-R_1^2=r_{in}^2(\hat{t}(l))-\hat{r}_{R_1}^2(l)
\end{equation}
where $\hat{r}_{R_1}^2(l)$ is the function $\hat{r}$ for this specific characteristic curve.
This curve divides the body into two regions.
As can be seen graphically from figure \ref{fig:outer_acc_charac}, the characteristic curves to the left of this specific curve (green region) originate from the line $t=0$; therefore, this part of the body is governed by the initial conditions of the body. 
On the other hand, the characteristic curves to the right of this specific curve (red region) originate from the outer boundary; therefore, this part of the body is governed by the boundary condition on the outer surface.
Consequently, $\bfF_e$ and $\bfsigma$ behave differently in each of these regions, and they need to be treated separately. 

\subsubsection{The region that was initially in the body}
    
Considering the region that formed part of the initial body, the characteristic curves from \eqref{eqn:outeracc_char_in} can be used to substitute $C$ in \eqref{eqn:Frr_general_outer} and \eqref{eqn:Ftt_general_outer} to find the general forms of $F_{e_{rr}}$ and $F_{e_{\theta\theta}}$:
\begin{align}
    \label{eqn:Frr_general_outer1}
   & \hat{r} F_{e_{rr}}(\hat{r},\hat{t}) = g_1(C_1) =
   g_1(r_{in}^2(\hat{t})-\hat{r}^2)
   \\
   \label{eqn:Ftt_general_outer1}
   & \frac{F_{e_{\theta\theta}}(\hat{r},\hat{t})}{\hat{r}} =
   g_2(C_1)=g_2(r_{in}^2(\hat{t})-\hat{r}^2)
\end{align}
Using the initial condition of $F_{e_{rr}}=F_{e_{\theta\theta}}=1$ at $t=0$ to compute $g_1$ and $g_2$, and substituting $\hat{r}$ and $\hat{t}$ with $r$ and $t$, we can find $\bfF_e$ in this region: 
\begin{equation}
    \label{eqn:def_grad_outer_1}
    \bfF_e=\frac{\sqrt{R_0^2-r_{in}^2(t)+r^2}}{r} \hat{\bfe}_r \otimes \hat{\bfe}_r +
    \frac{r}{\sqrt{R_0^2-r_{in}^2(t)+r^2}} \hat{\bfe}_\theta \otimes \hat{\bfe}_\theta
\end{equation}
Using \eqref{eqn:outeracc_char_in}, for a particle that was initially ($l=0$) at $R$, the spatial location at any later time $t$ is:
\begin{equation}
    \label{eqn:outer_r_R}
   r^2=R^2+r_{in}^2(t)-R_0^2
   \quad \Rightarrow \quad
   \sqrt{R_0^2-r_{in}^2(t)+r^2}=R
\end{equation}
So, we can write $\bfF_e$ as:
\begin{equation*}
    \bfF_e=\frac{R}{r} \hat{\bfe}_r \otimes \hat{\bfe}_r +
    \frac{r}{R} \hat{\bfe}_\theta \otimes \hat{\bfe}_\theta
\end{equation*}
We notice that this agrees with \cite{sozio2017nonlinear}, because the initial body is stress-free, and hence the initial configuration is also the reference configuration for this part of the body.

Substituting $\bfF_e$ from \eqref{eqn:def_grad_outer_1} into \eqref{eqn:const-law-outer}, the Cauchy stress tensor in this region is:
\begin{equation}
\label{eqn:outer_sig_region1}
\begin{split}
    \bfsigma=
    &\left[-p+2\parderiv{W}{I_2}\left( 1-\frac{r^2}{R_0^2-r_{in}^2(t)+r^2}\right)+2\parderiv{W}{I_1}\left(\frac{R_0^2-r_{in}^2(t)+r^2}{r^2} \right)\right]\hat{\bfe}_r \otimes \hat{\bfe}_r
    \\
    & + \left[-p+2\parderiv{W}{I_2}\left( 1-\frac{R_0^2-r_{in}^2(t)+r^2}{r^2}\right)+2\parderiv{W}{I_1}\left(\frac{r^2}{R_0^2-r_{in}^2(t)+r^2} \right)\right]\hat{\bfe}_\theta \otimes \hat{\bfe}_\theta
\end{split}
\end{equation}
To compute the unknown pressure field $p$, we use \eqref{eqn:mom-outer} and the appropriate boundary condition \eqref{eqn:bc-inner-mom-outer}:
\begin{equation}
\begin{split}
     &
     \parderiv{\sigma_{rr}}{r}+\frac{1}{r}(\sigma_{rr}-\sigma_{\theta\theta}) = 0 \Rightarrow
     \\
     &\parderiv{\sigma_{rr}}{r}=\frac{2}{r} \left(\parderiv{W}{I_1}+\parderiv{W}{I_2} \right) \left(\frac{r^4-\left(R_0^2-r_{in}^2(t)+r^2\right)^2}{\left(R_0^2-r_{in}^2(t)+r^2\right)r^2} \right)
     \\
     & \qquad 
     = \frac{2r}{R_0^2-r_{in}^2(t)+r^2} \left(\parderiv{W}{I_1}+\parderiv{W}{I_2} \right) \left(1 - \frac{(R_0^2-r_{in}^2(t)+r^2)^2}{r^4} \right)
\end{split}
\end{equation}
As $\sigma_{rr}$ is a function of $r$ and $t$, we have:
\begin{equation*}
    \sigma_{rr}=\int_{r_{in}(t)}^{r} \left[\frac{2r}{R_0^2-r_{in}^2(t)+r^2} \left(\parderiv{W}{I_1}+\parderiv{W}{I_2} \right) \left(1 - \frac{(R_0^2-r_{in}^2(t)+r^2)^2}{r^4} \right) \right] \dm r + h(t)
\end{equation*}
Using the boundary condition \eqref{eqn:bc-inner-mom-outer}, $h(t)$ is equal to $-p_i(t)$, thus:
\begin{equation}
\label{eqn:sig_rr_out1_1}
    \sigma_{rr}=-p_i(t) + \int_{r_{in}(t)}^{r} \left[\frac{2r}{R_0^2-r_{in}^2(t)+r^2} \left(\parderiv{W}{I_1}+\parderiv{W}{I_2} \right) \left(1 - \frac{(R_0^2-r_{in}^2(t)+r^2)^2}{r^4} \right) \right] \dm r
\end{equation}
And,
\begin{equation}
\label{eqn:outer_p}
\begin{split}
    p(r,t)=
    &
    p_i(t) -\int_{r_{in}(t)}^{r} \left[\frac{r}{R_0^2-r_{in}^2(t)+r^2} \left(\parderiv{W}{I_1}+\parderiv{W}{I_2} \right) \left(1 - \frac{(R_0^2-r_{in}^2(t)+r^2)^2}{r^4} \right) \right] \dm r
    \\
    & 
    + 2\parderiv{W}{I_2}\left( 1-\frac{R_0^2-r_{in}^2(t)+r^2}{r^2}\right)+2\parderiv{W}{I_1}\left(\frac{r^2}{R_0^2-r_{in}^2(t)+r^2} \right)        
\end{split}
\end{equation}
Using $p$ from \eqref{eqn:outer_p} in \eqref{eqn:outer_sig_region1} gives the complete Cauchy stress tensor in this region of the body.

To compare this result with \cite{sozio2017nonlinear}, we need to write \eqref{eqn:sig_rr_out1_1} based on the reference configuration and in Lagrangian form. Using \eqref{eqn:outer_r_R}, we have $r \dm r = R \dm R$, and substituting this relation in \eqref{eqn:sig_rr_out1_1}, we can write $\sigma_{rr}$
in the following form that agrees with \cite{sozio2017nonlinear}.
\begin{equation*}
    \sigma_{rr}=-p_i(t) + \int_{R_1}^{R} \left[\frac{2}{R} \left(\parderiv{W}{I_1}+\parderiv{W}{I_2} \right) \left(1 - \frac{R^4}{r^4} \right) \right] \dm R
\end{equation*}

\subsubsection{The region that is added to the body via growth process}
    
In the region that is governed by the outer boundary conditions, it is more convenient to use \eqref{eqn:outeracc_char_out} to compute the characteristic curves as it involves the outer part of the body, in contrast to \eqref{eqn:outeracc_char_in}.
So, for $C$ we use $C_2$ in \eqref{eqn:Frr_general_outer} and \eqref{eqn:Ftt_general_outer}.
Hence, the general forms of $F_{e_{rr}}$ and $F_{e_{\theta \theta}}$ are:
\begin{align}
   \label{eqn:Frr_general_outer2}
   & \hat{r} F_{e_{rr}}(\hat{r},\hat{t}) = g_1(C_2) =
   g_1\left(r_{out}^2(\hat{t})-\hat{r}^2-2\int_0^{\hat{t}} u_g(t) r_{out}(t) \dm t \right)
   \\
   \label{eqn:Ftt_general_outer2}
   & \frac{F_{e_{\theta\theta}}(\hat{r},\hat{t})}{\hat{r}} =
   g_2(C_2)=g_2\left(r_{out}^2(\hat{t})-\hat{r}^2-2\int_0^{\hat{t}}u_g(t) r_{out}(t) \dm t \right)
\end{align}
Using the boundary condition from \eqref{eqn:Fdot_bc_ic_outer} at $r=r_{out}(t)$, we have:
\begin{align}
    &
    \label{eqn:bc-rr-region2-outer}
    r_{out}(\hat{t})=g_1 \left(-2\int_0^{\hat{t}} u_g(t) r_{out}(t) \dm t \right)
    \\
    &
    \label{eqn:bc-tt-region2-outer}
    \frac{1}{r_{out}(\hat{t})}=g_2 \left(-2\int_0^{\hat{t}} u_g(t) r_{out}(t) \dm t \right)
\end{align}
To compute $F_{e_{rr}}$ and $F_{e_{\theta \theta}}$, we need to know the functions $g_1$ and $g_2$ when their arguments are $r_{out}^2(t)-\hat{r}^2-2\int_0^{t} u_g(t) r_{out}(t) \dm t$.
Using \eqref{eqn:outeracc_char_out} for a particle that is initially at the outer radius of $R_1$, substituting $(r=r_{out}(t=0)=R_1, t=0)$, the constant $C_2$ for this characteristic curve is computed to be zero.
Hence, another form of $\hat{r}_{R_1}(l)$ from \eqref{eqn:outer_rin_r(R_0)} is as follows: \begin{equation}
    \label{eqn:outer_r(R_0) out}
     \hat{r}^2_{R_1}(l) =r_{out}^2(\hat{t})-2\int_0^{\hat{t}} u_g(t) r_{out}(t) \dm t
\end{equation}
The schematic sketch of $\hat{r}_{R_1}(l)$ is the bold dashed line, and $r_{out}(t)$ is the bold dotted dashed line in figure \ref{fig:outer_acc_charac}.

Further, for a particle that is at $r=r_{out}(\tau)$ at $t=\tau$, substituting $(r=r_{out}(t=\tau) , t=\tau)$, the constant $C_2$ is $-2 \int_0^{\tau} u_g(t) r_{out}(t) \dm t$.
Thus, the pathline of a particle that attaches to the body at the outer radius $r_{out}(\tau)$ (at time $t=\tau$) at any time after the attachment time $t>\tau$ is:
\begin{equation}
    \label{eqn:outer_r()}
    r_{out}^2(\hat{t})- 2\int_0^{\hat{t}} u_g(t) r_{out}(t) \dm t = {\hat{r}}^2 - 2\int_0^{\tau} u_g(t) r_{out}(t) \dm t
\end{equation}
Using this equation and boundary conditions \eqref{eqn:bc-rr-region2-outer} and \eqref{eqn:bc-tt-region2-outer}, we get:
\begin{align}
    &
    r_{out}^2(\hat{t})-\hat{r}^2-2\int_0^{\hat{t}}u_g(t) r_{out}(t) \dm t= -2\int_0^{\tau}u_g(t) r_{out}(t) \dm t \Rightarrow 
    \\
    &
    \hat r F_{e_{rr}}(\hat{r}, \hat{t}) = g_1 \left(r_{out}^2(\hat{t})-\hat{r}^2-2\int_0^{\hat{t}}u_g(t) r_{out}(t) \dm t \right) = g_1 \left(-2\int_0^{\tau}u_g(t) r_{out}(t) \dm t \right) = r_{out}(\tau)
    \\
    &
    \frac{F_{e_{\theta \theta}}(\hat{r}, \hat{t})}{\hat r} = g_2 \left(r_{out}^2(\hat{t})-\hat{r}^2-2\int_0^{\hat{t}}u_g(t) r_{out}(t) \dm t \right) = g_2 \left(-2\int_0^{\tau}u_g(t) r_{out}(t) \dm t \right) = \frac{1}{r_{out}(\tau)}
\end{align}
So, the deformation gradient in this region is:
 \begin{equation}
\label{eqn:outer_def_grad_2}
    \bfF_e=\frac{r_{out}(\tau (\hat{r}, \hat{t}))}{r} \hat{\bfe}_r \otimes \hat{\bfe}_r +
    \frac{r}{r_{out}(\tau(\hat{r}, \hat{t}))} \hat{\bfe}_\theta \otimes \hat{\bfe}_\theta
\end{equation}
Now, we need to compute $\tau$ for a given spatial location $r$ at given time $t$.
Substituting \eqref{eqn:outer_r()} into \eqref{eqn:outer_r(R_0) out}, the following relation implicitly provides an equation to compute $\tau$ for a particle that is at the spatial radius $r$ at time $t$ with $\hat{r}_{R_1}(t) < r < r_{out}(t)$:
\begin{equation}
\label{eqn:outer_r_R_0_tau}
    \hat{r}^2(l) = \hat{r}^2_{R_1} + 2\int_0^{\tau} u_g(t) r_{out}(t) \dm t
\end{equation}
The (as yet unknown) functions $r_{out}(\hat{t})$ and $\hat{r}_{R_1}$ are required apply this equation to compute $\tau$. 
The functions $r_{out}(\hat{t})$ and $\hat{r}_{R_1}$ are related to each other by \eqref{eqn:outer_r(R_0) out}, so indeed we have only one unknown function. 
For now, we will continue with these unknown functions, and they will be computed using the traction boundary conditions below.

Although we want to compute $\tau$ for a known spatial location $r$ at known time $t$ because of the Eulerian approach, \cite{sozio2017nonlinear} used a similar equation to compute the unknown location $r$ at known time $t$ of a particle that has reference position $R$ in the Lagrangian approach.
In that approach, one needs to also know the relation between $\tau$ and $R$ as the attachment time affects the current location of the body, and the integral in the last term of the above equation is written with respect to $\dm R$ in \cite{sozio2017nonlinear}. 
For simplicity, \cite{sozio2017nonlinear} assume that $u_g(t)$ is constant and equal to $u_0$, and then $\tau = \frac{R-R_1}{u_0}$.
In general, the relation between $\tau$ to $R$ is more complicated. 
The Eulerian approach does not refer to the reference configuration explicitly, and hence the choice of $u_g(t)$ does not affect the formulation of the problem.

We notice that we could alternatively find $\bfF_e$ from \eqref{eqn:outer_def_grad_2} using the following argument.
The boundary condition for this region is that the particles are added to the body at the outer surface with radius of $r_{out}(\hat{t})$ without experiencing any deformation, i.e. $\bfF_e=\bfI$. 
In the other words, $F_{e_{rr}} = F_{e_{\theta \theta}} = 1$ at the attachment radius of the particle $r=r_{attach}(\hat{r},\hat{t})$; this is different for each characteristic curve.
Also, we know that along pathlines, $\hat{r} F_{e_{rr}}$ and $\frac{F_{e_{\theta\theta}}}{\hat{r}}$ are constants and equal to $g_1(C_2)$ and $g_2(C_2)$, respectively. 
Using the boundary condition, $r_{attach}(\hat{r},\hat{t})=g_1(C_2)=\frac{1}{g_2(C_2)}$, and we can then write $\bfF_e$ in this region as:

\begin{equation*}
    \bfF_e=\frac{r_{attach}(\hat{r},\hat{t})}{r} \hat{\bfe}_r \otimes \hat{\bfe}_r +
    \frac{r}{r_{attach}(\hat{r},\hat{t})} \hat{\bfe}_\theta \otimes \hat{\bfe}_\theta
\end{equation*}
And obviously the attachment radius of each particle is the outer radius of the body at the time of attachment $\tau$, so $r_{attach}(\hat{r},\hat{t})=r_{out}(\tau)$; \eqref{eqn:outer_r_R_0_tau} gives us an equation to compute $\tau$. 
    
Now that $\bfF_e$ is known \eqref{eqn:outer_def_grad_2}, we can compute the Cauchy stress using \eqref{eqn:const-law-outer}:
\begin{equation}
\label{eqn:outer_sig_region2}
\begin{split}
    \bfsigma=
    &\left[-p+2\parderiv{W}{I_2}\left( 1-\frac{r^2}{r_{out}^2(\tau)}\right)+2\parderiv{W}{I_1}\left(\frac{r_{out}^2(\tau)}{r^2} \right)\right]\hat{\bfe}_r \otimes \hat{\bfe}_r
    \\
    & + \left[-p+2\parderiv{W}{I_2}\left( 1-\frac{r_{out}^2(\tau)}{r^2}\right)+2\parderiv{W}{I_1}\left(\frac{r^2}{r_{out}^2(\tau)} \right)\right]\hat{\bfe}_\theta \otimes \hat{\bfe}_\theta        
\end{split}
\end{equation}
To find the unknown pressure field $p$, we use \eqref{eqn:mom-outer} to compute $\sigma_{rr}$:
\begin{equation*}
    \sigma_{rr}
    =
    \int_{\hat{r}_{R_1}}^{r} \left[\frac{1}{r} (\sigma_{\theta \theta} - \sigma_{rr}) \right] \dm r
    =
    \int_{\hat{r}_{R_1}}^{r} \left[\frac{2r}{r_{out}^2(\tau)} \left(\parderiv{W}{I_1}+\parderiv{W}{I_2} \right) \left(1 - \frac{r_{out}^4(\tau)}{r^4} \right) \right] \dm r + h(t)
\end{equation*}
Using traction continuity at the interface between the initial and added material (located at $\hat{r}_{R_1}$), we have for $\sigma_{rr}$ and the function $h(t)$:
\begin{align*}
    &
    h(t) =
    -p_i(t) +
    \int_{r_{in}(t)}^{r(R_1,t)} \left[\frac{2r}{R_0^2-r_{in}^2(t)+r^2} \left(\parderiv{W}{I_1}+\parderiv{W}{I_2} \right) \left(1 - \frac{(R_0^2-r_{in}^2(t)+r^2)^2}{r^4} \right) \right] \dm r
    \\
    \sigma_{rr}=
    &
    -p_i(t) +
    \int_{r_{in}(t)}^{r(R_1,t)} \left[\frac{2r}{R_0^2-r_{in}^2(t)+r^2} \left(\parderiv{W}{I_1}+\parderiv{W}{I_2} \right) \left(1 - \frac{(R_0^2-r_{in}^2(t)+r^2)^2}{r^4} \right) \right] \dm r+
    \\
    &
    \int_{r(R_1,t)}^{r} \left[\frac{2r}{r_{out}^2(\tau)} \left(\parderiv{W}{I_1}+\parderiv{W}{I_2} \right) \left(1 - \frac{r_{out}^4(\tau)}{r^4} \right) \right] \dm r  
\end{align*}
We note that $r_{in}(t)$, $r_{out}(t)$ and $\hat{r}_{R_1}$ are as yet unknown functions. 
Based on \eqref{eqn:in_out_vel_relation_outeracc} and \eqref{eqn:outer_rin_r(R_0)}, these functions are not independent, and we can write both $\hat{r}_{R_1}$ and $r_{out}(t)$ in terms of $r_{in}(t)$.
Another boundary condition is that the outer surface $r=r_{out}(t)$ is stress-free. 
We therefore get the following equation for computing these three (non-independent) unknown functions:
\begin{equation}
\begin{split}
    p_i(t) =
    & \int_{r_{in}(t)}^{r(R_1,t)} \left[\frac{2r}{R_0^2-r_{in}^2(t)+r^2} \left(\parderiv{W}{I_1}+\parderiv{W}{I_2} \right) \left(1 - \frac{(R_0^2-r_{in}^2(t)+r^2)^2}{r^4} \right) \right] \dm r
    \\
    & + \int_{r(R_1,t)}^{r_{out}(t)} \left[\frac{2r}{r_{out}^2(\tau)} \left(\parderiv{W}{I_1}+\parderiv{W}{I_2} \right) \left(1 - \frac{r_{out}^4(\tau)}{r^4} \right) \right] \dm r
\end{split}
\end{equation}
This is in agreement with \cite{sozio2017nonlinear}.

\subsection{Computing the deformation gradient using inverse of the motion}

An alternative Eulerian approach to find $\bfF$ is to work with the inverse of the motion $\bfxi=\xi_r \hat{\bfe}_r+\xi_{\theta} \hat{\bfe}_{\theta}$ and \eqref{eqn:Kamrin} provided in \cite{kamrin2009eulerian}; once we have $\bfxi$, we can compute $\bfF$ using $\bfF=(\nabla\bfxi)^{-1}$. And we note that as $\bfFr = \bfI$, $\bfF_e = \bfF$. 
\begin{equation}
\label{eqn:Kamrin}
    \parderiv{\bfchi^{-1}}{t} + (\bfv\cdot\nabla) \bfchi^{-1} = 0
\end{equation}
From symmetry, we have $\xi_{\theta} = \const$.
For $\xi_r$, we have:
\begin{equation}
    \parderiv{\xi_r}{t} + \dot{r}_{in} \frac{r_{in}}{r} \parderiv{\xi_r}{r} = 0
\end{equation}
which is the $r-$component of \eqref{eqn:Kamrin}, after substituting the velocity field from \eqref{eqn:vel_outer}.

Using the method of characteristics, the characteristic curves are the same as those from \eqref{eqn:outeracc_char_in} and \eqref{eqn:outeracc_char_out}. 
However, as the right side of the above equation is zero, it means that $\xi_r$ is constant along the characteristics, giving:
\begin{equation}
    \label{eqn:xi_general_outer}
    \xi_r=g(C)=
    g_1 \left(r_{in}^2(\hat{t})-\hat{r}^2 \right)=
    g_2 \left(r_{out}^2(\hat{t}) - \hat{r}^2 - 2 \int_{0}^{\hat{t}} u_g(t) r_{out}(t) \dm t \right)
\end{equation}
We can then compute $\bfF$ in the initial region:
\begin{equation}
    \bfF=(\nabla\bfxi)^{-1}
    =
    \begin{bmatrix}
        \parderiv{\xi_r}{r} & 0 \\
        0 & \frac{\xi_r}{r}
    \end{bmatrix} ^{-1} 
    =
    \begin{bmatrix}
        \frac{1}{-2 r g_1'(m)} & 0 \\
        0 & \frac{r}{g_1(m)}
    \end{bmatrix}_{m=r_{in}^2(t)-r^2}
\end{equation}
This region is influenced by the initial condition, which is $\bfF=\bfI$.
Also, the intersection between the line $t=0$ and a characteristic curve that passes through $\hat{r}$ at time $\hat{t}$ is $R$ (the initial location of the particle), so:
\begin{equation}
    g_1(R_0^2-R^2)=R \Rightarrow g_1(m)=\sqrt{R_0^2 - m} 
    \quad \Rightarrow \quad 
    g_1'(m)= \frac{-1}{2 \sqrt{R_0^2 - m}}
\end{equation}
Therefore, the explicit expression for $\bfF$ is:
\begin{equation}
    \bfF=
    \begin{bmatrix}
    \frac{\sqrt{R_0^2-r_{in}^2(t)+r^2}}{r} & 0 \\
    0 & \frac{r}{\sqrt{R_0^2-r_{in}^2(t)+r^2}}
    \end{bmatrix}
\end{equation}
which is, as expected, identical to \eqref{eqn:def_grad_outer_1}.

In the region governed by the boundary condition, we have from \eqref{eqn:xi_general_outer} that $\xi_r=g_2(C_2)$.
Combined with the boundary condition $\bfF=\bfI$ at $r=r_{out}(t)$, we have:
\begin{equation}
    \bfF=
    \begin{bmatrix}
    \frac{1}{-2 r g_2'(m)} & 0 \\
    0 & \frac{r}{g_2(m)}
    \end{bmatrix}_{m=\left(r_{out}^2(\hat{t}(l)) - \hat{r}^2(l,f) - 2 \int_{0}^{\hat{t}(l)} u_g(t) r_{out}(t) \dm t\right)}
\end{equation}
At $r=r_{out}(\hat{t}),$ we have $F_{\theta\theta}=1$, so $r_{out}(\hat{t})=g_2\left(- 2 \int_{0}^{\hat{t}} u_g(t) r_{out}(t) \dm t\right)$. 
Using \eqref{eqn:outer_r()}:
\begin{equation*}
    g_2(m)|_{m=\left(r_{out}^2(\hat{t}) - \hat{r}^2 - 2 \int_{0}^{\hat{t}} u_g(t) r_{out}(t) \dm t\right)} = 
    g_2\left(- 2 \int_{0}^{\tau} u_g(t) r_{out}(t) \dm t\right) = 
    r_{out}(\tau (\hat{r}, \hat{t}))
\end{equation*}
Now, we need $g_2'(m)$ to define $\bfF$ completely; however, the general form of $g_2(m)$ is not readily available to differentiate.
We only know that $g_2(m)=r_{out}(\tau)$ at $m=r_{out}^2(\hat{t}(l)) - \hat{r}^2 - 2 \int_{0}^{\hat{t}(l)} u_g(t) r_{out}(t) \dm t$. 
Therefore, we consider $g'(m)$ as an unknown function, and apply the boundary condition $F_{rr}=1$ at $r=r_{out}$ to compute $g_2'(m)$ directly:
\begin{align*}
    &-2 r g_2'\left( - 2 \int_{0}^{\hat{t}(l)} u_g(t) r_{out}(t) \dm t\right)=1 \Rightarrow
    g_2'\left( - 2 \int_{0}^{\hat{t}} u_g(t) r_{out}(t) \dm t\right)=
    \frac{1}{-2 r_{out}(\hat{t})} \Rightarrow \\
    &g'_2\left(r_{out}^2(\hat{t}(l)) - \hat{r}^2 - 2 \int_{0}^{\hat{t}(l)} u_g(t) r_{out}(t) \dm t\right)=
    g'_2\left(- 2 \int_{0}^{\tau} u_g(t) r_{out}(t) \dm t\right)=
    \frac{1}{-2 r_{out}(\tau)} \Rightarrow\\
    & F_{rr}=\frac{r_{out}(\tau)}{r}
\end{align*}
Therefore, $\bfF$ in the region governed by the boundary conditioned region is:
\begin{equation}
    \bfF=
    \begin{bmatrix}
    \frac{r_{out}(\tau)}{r} & 0 \\
    0 & \frac{r}{r_{out}(\tau)}
    \end{bmatrix}
\end{equation}
which is identical to \eqref{eqn:outer_def_grad_2}. 

Instead of using the initial condition on $\bfF$, we could use that the region governed by the initial condition is initially stress-free to set $\xi_r=R$, and use this to compute the function $g_1(C)$ in the initially existing body:
\begin{equation*}
    g_1(R_0^2-R^2)=R \Rightarrow \xi_r = g_1(C_2)=\sqrt{R_0^2 - r_{in}^2(\hat{t}) - \hat{r}^2}    
\end{equation*}
However, in the region governed by the boundary condition, converting the boundary condition on $\bfF$ to a boundary condition on $\xi_r$ requires the initial location of the added particles, and hence has the same difficulty as computing the reference configuration.
This is the reason that we compute $\bfF$ using \eqref{eqn:def-grad-transport} rather than through computing the inverse of the motion using \eqref{eqn:Kamrin}.



\section{Discussion}
\label{sec:conclusion}

In this paper, we model the surface growth of solid bodies using an Eulerian approach.
An important advantage of this approach is that we do not need to explicitly compute the time-evolving reference configuration.
On the other hand, to obtain the stress, we need to solve for the deformation gradient which is not straightforward in the Eulerian setting.
This is addressed by introducing the elastic deformation field and the associated transport equation that governs its evolution.
To solve this time-dependent transport equation, we use the method of characteristics that enables us to find closed-form solutions in two examples with spherical and cylindrical symmetry with free boundaries.

The primary motivation for formulating the Eulerian approach is to enable future work on numerical solutions.
A potential advantage of Eulerian methods that have been developed in the fluid-structure interaction (FSI) literature (e.g. \cite{sugiyama2011full}) is that we can solve such problems on a fixed mesh for an evolving body, rather than an evolving reference configuration and deforming mesh that needs to be updated or re-meshed frequently.
A numerical scheme would enable the careful study of various interesting systems identified recently, e.g. \cite{zurlo2017printing, zurlo2018inelastic, truskinovsky2019nonlinear} that focus on the problem of near-net-shape manufacturing in the presence of incompatibility and residual stress, as well as \cite{swain2018biological} that discusses incompatibility in both surface and bulk growth.

    An important challenge in the Eulerian formulation is to know the location of the boundaries to be able to apply the boundary conditions.
    We refer to the brief discussion in Section 3.3.1 of \cite{naghibzadeh2021surface}, as well as \cite{kamrin2012reference,kamrin2009eulerian}, and notice that the velocity field can be used to locate the boundary of the body
    for the parts of the boundary without growth, we follow the characteristic lines of particles located on the boundary; for the parts with growth, we use \eqref{eqn:bc-balance-mass} together with the pathlines to determine the location of the growing boundaries.

Briefly, we also highlight that bulk growth has received much attention, e.g. using growth tensors \cite{rodriguez1994stress, chenchiah2014energy, reina2017incompressible,garikipati2009kinematics}, geometrical approaches \cite{yavari2010geometric}, and mixture theory \cite{humphrey2002constrained}.
A key difference between bulk and surface growth is that the former can typically be studied by assuming that the set of the material points is unchanged, and that growth occurs through a change in the \textit{referential} density due to a distributed mass source.

\begin{acknowledgments}
    We thank Rohan Abeyaratne, Timothy Breitzman, Tal Cohen, and Tony Rollett for useful discussions.
    This work was supported by the National Science Foundation (CMMI MOMS 1635407, DMREF 2118945, DMS 1729478, DMS 2012259, DMS 2108784), Army Research Office (MURI W911NF-19-1-0245), Office of Naval Research (N00014-18-1-2528), and Binational Science Foundation (2018183).
\end{acknowledgments}

\bibliographystyle{alpha}
\bibliography{growth-refs}

\end{document}